\documentclass[12pt]{amsart}
\usepackage{amsfonts}
\usepackage{graphicx}
\usepackage{amsmath, amscd, amssymb, yhmath}
\usepackage{latexsym}
\usepackage{mathrsfs}
\usepackage{xy}
\input xy
\xyoption{all} 

\def \R{\mathbb{R}}

\def \C{\mathbb{C}}

\def \J{\mathcal{J}}

\def \K{\mathcal{K}}
\def \mo{\mathop{\mathrm{MO}}}
\def \MO{\mathbf{MO}}
\def \MSO{\mathbf{MSO}}
\def \BDiff{\mathop{\mathrm{BDiff}}}
\def \Diff{\mathop{\mathrm{Diff}}}
\def \BO{\mathop{\mathrm{BO}}}
\def \Ob{\mathop{\mathrm{Ob}}}
\def \Mor{\mathop{\mathrm{Mor}}}

\def \id{\mathrm{id}\,}
\newtheorem{theorem}{Theorem}[section]
\newtheorem{lemma}{Lemma}[section]

\newtheorem{corollary}{Corollary}

\theoremstyle{remark}
\newtheorem{remark}[theorem]{Remark}
\newtheorem{convention}[theorem]{Convention}

\theoremstyle{definition}
\newtheorem{definition}{Definition}

\newtheorem{example}{Example}

\begin{document}

\title{Singular cobordism categories}
\author{Rustam Sadykov}
\address{Mathematics Department\\ University of Toronto\\ Canada}
\email{rstsdk@gmail.com}
\thanks{The author has been supported by
Postdoctoral Fellowships at Max Planck Institute, Germany, and at the University of Toronto, Canada}
\date{\today}

\begin{abstract} Recently Galatius, Madsen, Tillmann and Weiss identified the homotopy type of the classifying space of the cobordism category of embedded $d$-dimensional manifolds \cite{GMTW} for each positive integer $d$. Their result lead to a new proof of the generalized standard Mumford conjecture. We extend the main theorem of \cite{GMTW} to the case of cobordism categories of embedded $d$-dimensional manifolds with prescribed singularities, and explain the relation of singular cobordism categories to the bordism version of the Gromov h-principle.   
\end{abstract}

\maketitle

\setcounter{tocdepth}{1}
\tableofcontents

\addcontentsline{toc}{part}{Introduction}

\section{Introduction}

Recently Galatius, Madsen, Tillmann and Weiss identified the homotopy type of the classifying space of the cobordism category of embedded $d$-dimensional manifolds \cite{GMTW} for each positive integer $d$. We extend their theorem to the case of cobordism categories of embedded $d$-dimensional manifolds with prescribed singularities. 
\begin{figure}[ht]
\includegraphics[draft=false, width=80mm]{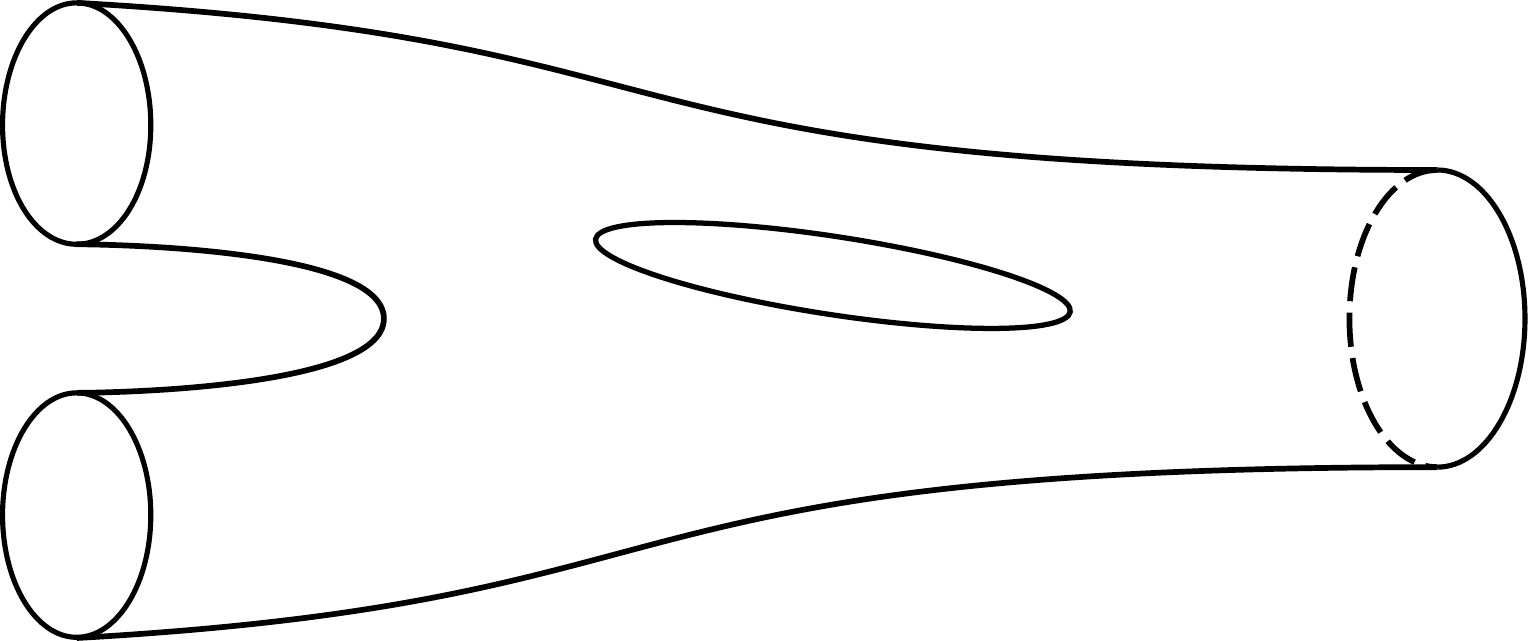}
\caption{A morphism $W$ between $M_1$ and $M_2$ in $\mathcal{C}_d$ is an (embedded) smooth cobordism between $M_1$ and $M_2$.}
\end{figure}

The \emph{cobordism category $\mathcal{C}_d$ of smooth manifolds} of dimension $d$ is a generalization of the category of conformal surfaces introduced by Segal. An object in $\mathcal{C}_d$ is a closed smooth manifold of dimension $d-1$ embedded into $\{a\}\times \R^{\infty+d-1}$ for some real number $a$, while a non-identity morphism in $\mathcal{C}_d$ from a manifold $M_1$ in $\{a_1\}\times \R^{\infty+d-1}$ to a manifold $M_2$ in $\{a_2\}\times \R^{\infty+d-1}$ exists only if $a_1<a_2$, in which case it is given by a compact smooth submanifold   
\[
W\subset [a_0, a_1]\times \R^{\infty+d-1}
\] 
of dimension $d$ transversally intersecting the walls $\{a_0, a_1\}\times \R^{\infty+d-1}$ in $M_1\sqcup M_2$ (see section~\ref{se:4}). The composition in the category $\mathcal{C}_d$ is defined by taking the union of submanifolds $W$. 

It turned out \cite{GMTW} that the loop space $\Omega B\mathcal{C}_2$ of the classifying space of $\mathcal{C}_2$ is rationally homology equivalent to the stable moduli space of Riemann surfaces, which allowed the authors of \cite{GMTW} to give a new proof of the standard Mumford conjecture, which originally was solved by Madsen and Weiss in \cite{MW} (see also \cite{EGM}, \cite{GR}). According to a celebrated work of Pierre Deligne and David Mumford, each moduli space of Riemann surfaces admits a  compactification by a moduli space, called the {\it Deligne-Mumford compactification}, of smooth Riemann surfaces as well as surfaces with so-called node singularities. The corresponding ``compactification" of cobordism categories of smooth manifolds motivates the central notion of the current paper, namely, the notion of cobordism categories $\mathcal{C}_{\J}$ of singular manifolds. 

\begin{figure}[ht]
\includegraphics[draft=false, width=80mm]{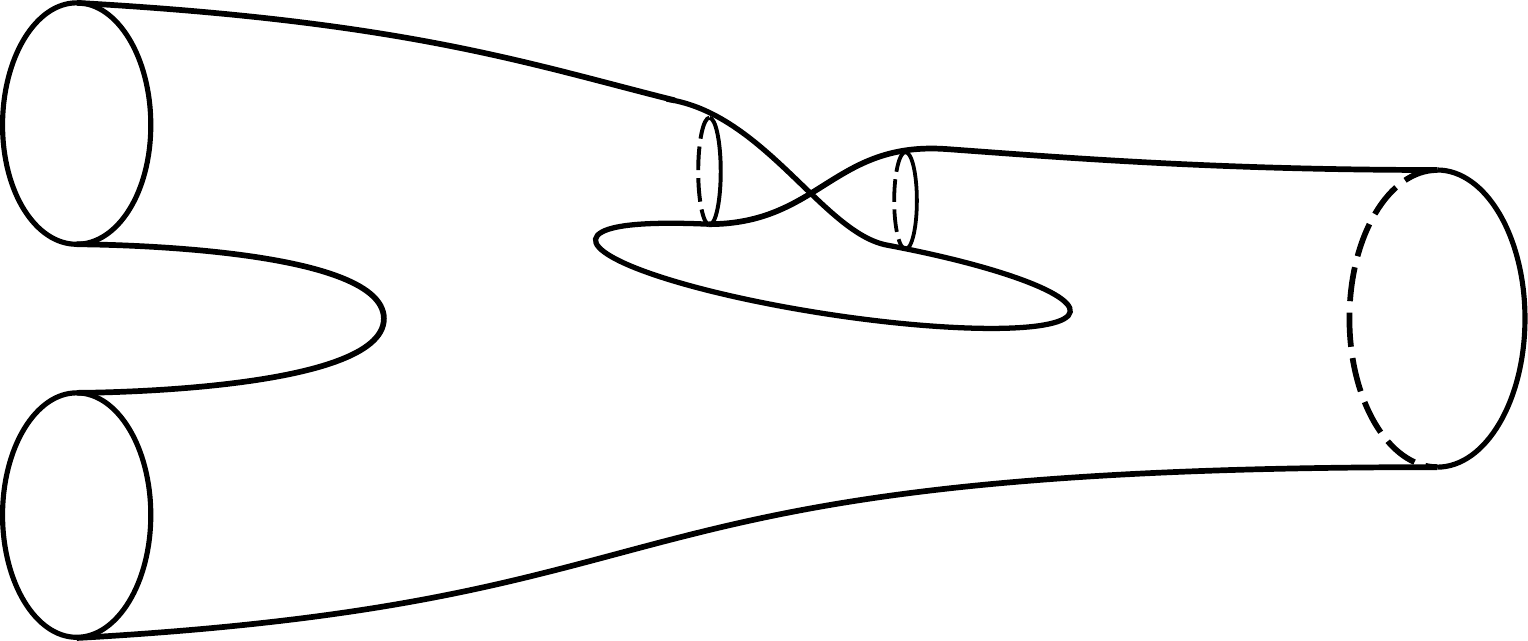}
\caption{A morphism in $\mathcal{C}_{\J}$; objects in $\mathcal{C}_{\J}$ are (embedded) smooth manifolds, while morphisms $W$ in $\mathcal{C}_{\J}$ are allowed to have singularities.}
\end{figure}

A \emph{singular manifold} of dimension $d$ is defined to be a fiber of a smooth map $M\to N$ of manifolds with $\dim M-\dim N=d$. For example, a singular manifold with Morse singularities is a fiber of a Morse function. 
Given a set  $\J$ of types of singularities of singular manifolds of dimension $d$, the {\it singular cobordism category $\mathcal{C}_{\J}$} is defined to be the category whose objects are the same as those in $\mathcal{C}_d$, i.e., embedded closed manifolds of dimension $d-1$, while the space of non-trivial morphisms in $\mathcal{C}_{\J}$ between embedded closed manifolds $M_1$ in $\{a_1\}\times \R^{\infty+d-1}$ and $M_2$ in $\{a_2\}\times \R^{\infty+d-1}$ consists of embedded compact singular manifolds $W$ in $[a_1, a_2]\times \R^{\infty+d-1}$ bounded by $\partial{W}=M_1\sqcup M_2$ with singularities of types in $\J$ (see section~\ref{s:9}). 

The main theorem of Galatius, Madsen, Tillmann and Weiss in \cite{GMTW} establishes a weak homotopy equivalence 
\begin{equation}\label{eq:i1}
   B\mathcal{C}_d \simeq \Omega^{\infty-1}\mathbf{B}_{\emptyset} 
\end{equation}
of the classifying space of $\mathcal{C}_d$ and a fairly simple infinite loop space $\Omega^{\infty}\mathbf{B}_{\emptyset}$. We extend this theorem to the case of cobordism categories of embedded manifolds with prescribed singularities under a fairly mild assumption that
the set of prescribed singularities $\J$ is \emph{open} and \emph{stable} (see section~\ref{sec:3}). In this case the counterpart of $\Omega^{\infty-1}\mathbf{B}_{\emptyset}$ is a fairly simple infinite loop space $\Omega^{\infty-1}\mathbf{B}_{\J}$ (see section~\ref{sec:3}). \\


\begin{theorem}~\label{th:1.4} Let $\J$ be an open stable set of singularities. Then there is a weak homotopy equivalence $B\mathcal{C}_{\J}\simeq \Omega^{\infty-1}\mathbf{B_{\J}}$.
\end{theorem}


Let $F_{g,1}$ denote an oriented surface of genus $g$ with one boundary component. Then $F_{g,1}\subset F_{g+1,1}$ and every orientation preserving diffeomorphism of $F_{g,1}$ pointwise trivial on the boundary $\partial F_{g,1}$ extends to an orientation preserving diffeomorphism of $F_{g+1,1}$. In particular there are inclusions $\Diff^+ F_{g,1}\subset \Diff^+ F_{g+1,1}$ of groups of orientation preserving diffeomorphisms pointwise trivial on the boundary. The colimit of the corresponding inclusions of classifying spaces  
\begin{equation}
    {\BDiff}^+ F_{0,1}\subset {\BDiff}^+ F_{1,1} \subset {\BDiff}^+ F_{2,1} \subset \cdots 
\end{equation}
is denoted by $\BDiff^+ F_{\infty, 1}$. Theorem~\ref{th:1.4} extends to the case of cobordism categories $\mathcal{C}_d^+$ and singular cobordism categories $\mathcal{C}_{\J}^+$ of \emph{oriented} manifolds; and in the case where $\J$ is empty, i.e., $\mathcal{C}^+_{\J}=\mathcal{C}^+_d$,  Galatius, Madsen, Tillmann and Weiss proved~\cite{GMTW} that there is a map
\begin{equation}\label{eq:1.5}
   \sqcup {\BDiff}^+ F_{\infty,1}\longrightarrow \Omega B\mathcal{C}_{2}^+,
\end{equation}
that induces an isomorphism of integral homology groups, which, in view of (\ref{eq:i1}), is equivalent to the generalized Mumford Conjecture. Let us now turn to the case of an arbitrary dimension $d\ge 0$ and a set $\J$ that in addition to hypotheses of Theorem~\ref{th:1.4} contains all Morse singularity types. In this case the counterpart of the homology equivalence (\ref{eq:1.5}) is a natural decomposition of the space $\Omega B\mathcal{C}^+_{\J}$.

\begin{theorem} \label{th:main}
The space $\Omega B\mathcal{C}^+_{\J}$ breaks into the union of subspaces
\begin{equation}\label{eq:1.6}
    \Omega B\mathcal{C}^+_{\J}=\cup {\BDiff}^+M_{\alpha}
\end{equation}
where the union ranges over the classifying spaces $\BDiff^+M_{\alpha}$ of groups $\Diff^+M_{\alpha}$ of orientation preserving diffeomorphisms of singular manifolds of dimension $d$. 
\end{theorem}

We emphasize that the components $\BDiff^+ M_{\alpha}$ of $\Omega B\mathcal{C}_{\J}$ are extremely complicated; these are not well understood already in the case where $M_{\alpha}$ is a smooth surface. On the other hand, Theorem~\ref{th:1.4} implies that the total space $\Omega B\mathcal{C}_{\J}$ is equivalent to the  space $\Omega^{\infty}\mathbf{B}_{\J}$, which is relatively simple. For example, the space $\Omega^{\infty}\mathbf{B}_{\emptyset}$ is constructed from the space $\BO_d$. We will see that the group $\mathop{\mathrm{O}}_d$ appears here as the symmetry group of a neighborhood of a point in a manifold of dimension $d$. In general, the space $\Omega^{\infty}\mathbf{B}_{\J}$ is build of the classifying spaces  
${\BDiff}\tau$ of symmetry groups of singularities $\tau\in\J$ of singular manifolds and the classifying space $\BO_d$ of the symmetry group of a neighborhood of a non-singular point in a smooth manifold of dimension $d$ (see Remark~\ref{symmetry}).  


\begin{remark} In the case where $\J$ consists of finitely many singularity types, the space $\Omega^{\infty}\mathbf{B}_{\J}$ can be constructed from spaces $\BDiff\tau$ as described in the paper \cite{EG} by Eliashberg and Galatius in the case  corresponding to the Deligne-Mumford compactification. The statement similar to that in Theorem~\ref{th:1.4} is known to be true for many types of cobordism categories. Furthermore, singular cobordism categories $\mathcal{C}_{\J}$ are closely related to cobordism categories of manifolds with tangential structures of Galatius, Madsen, Tillmann and Weiss \cite{GMTW} (see Remark~\ref{rem:1.7}), and geometric cobordism categories of Ayala~\cite{Aya}. 

On the other hand, versions of the equivalence (\ref{eq:1.5}) has been established only in those cases where there are stability theorems such as stability theorems by Harer~\cite{Ha}, Ivanov~\cite{Iv}, Wahl~\cite{Wa}, Galatius~\cite{Ga}, and Cohen-Madsen~\cite{CM} (see also a general approach by Galatius and Randal-Williams in \cite{GR}). We note that in all mentioned cases the dimension $d$ is $2$. Theorem~\ref{th:main} is true for an arbitrary dimension $d$.     
\end{remark}

\begin{remark}\label{rem:1.7} For a set $\J$ of singularity types as in the statement of Theorem~\ref{th:1.4}, there is a tangential structure $\theta$ on smooth manifolds such that $B\mathcal{C}_{\J}\simeq B\mathcal{C}_{\theta}$, where $\mathcal{C}_{\theta}$ is the cobordism category of manifolds with tangential structure $\theta$ \cite{GMTW}.  The natural ``tangential" structure $\theta$ on objects and morphisms of $\mathcal{C}_{\J}$, however, equips all objects and smooth morphisms with the trivial structure. Thus the obvious map $\Mor\mathcal{C}_{\J}\dashrightarrow \Mor\mathcal{C}_{d}$ is defined only on the subspace of smooth morphisms, while the obvious map $\Mor\mathcal{C}_{\theta}\dashrightarrow \Mor\mathcal{C}_{\J}$ is defined only on the subspace of morphisms with the trivial tangential structure.   
\end{remark}

We apply a technique developed in \cite{MW} and \cite{GMTW}, and we adopt the argument of Galatius, Madsen, Tillmann and Weiss \cite{GMTW} to the case of singular cobordism categories. A few essential modifications, however, are necessary. First, we introduce a suitable topology on the space of morphisms of categories $\mathcal{C}_{\J}$. This is done by a construction completely different from the one in \cite{GMTW}. Namely, we rely on the construction of Kazarian~\cite{Kaz1}, \cite{Kaz2} (see also \cite{SY}) in global singularity theory. Second, for a manifold bundle $f: M\to N$ and any smooth map $N'\to N$ there is a well-defined pullback manifold bundle $f': M'\to N'$. This allows the authors of \cite{GMTW} to introduce sheaves of manifold bundles on the category of smooth manifolds without boundary and smooth maps. For a general singular cobordism category such a pullback property is not available (see Remark~\ref{r:3.5}). Finally, the construction of the decomposition (\ref{eq:1.6}) is completely different from the construction of the map (\ref{eq:1.5}) in \cite{MW} and \cite{GMTW}. 

\subsection*{Acknowledgments}
The author is thankful to Hugo Chapdelaine for many comments on the preliminary version of the paper. The author is also thankful to Prof. Ib Madsen for an encouraging discussion of the results of the paper and hospitality at Copenhagen University.   

\section{Outline}

For the reader's convenience, in section~\ref{sec:3} we briefly review necessary notions from global singularity theory including the notion of a \emph{singularity of a smooth map}, and the notion of the \emph{type of a singularity}. For a set $\J$ of types of singularities of smooth maps, a \emph{$\J$-map} is a map with singularities of types in $\J$. A \emph{singular manifold} with $\J$-singularities is a fiber of a $\J$-map. In section~\ref{sec:3} we also define an infinite loop space $\Omega^{\infty}\mathbf{B}_{\J}$ and introduce sheaves $D_{\J}$ of singular manifolds. 



In section~\ref{s:6} we recall the Kazarian construction of the set $\Omega^{\infty}_{\J}\MO$. It is a subset of the infinite loop space $\Omega^{\infty+d}\MO$ defined so that a smooth proper map $f: M\to N$ of dimension $d$ is a $\J$-map if and only if the image of the associated Pontrjagin-Thom map $N\to \Omega^{\infty+d}\MO$ is in $\Omega^{\infty}_{\J}\MO$. 

\begin{example} In the case where $\J$ is empty, the set $\Omega^{\infty}_{\J}\MO$ consists of loops $S^{\infty+d}\to \MO$ transversal to the zero section $\BO\subset \MO$. A smooth proper map $f: M\to N$ is a $\J$-map if and only if the image of the associated Pontrjagin-Thom map is in $\Omega^{\infty+d}_{\emptyset}\MO$.   
\end{example}

In fact, in section~\ref{s:6} we define sets $\Omega^{\infty-k}_{\J}\MO$ for each non-negative integer $k$; and introduce a suitable topology on each of these spaces, including spaces $\Omega^{\infty}_{\J}\MO$. In section~\ref{s:7} we show that the space $\Omega^{\infty-1}\mathbf{B}_{\J}$ is weakly homotopy equivalent to $\Omega^{\infty-1}_{\J}\MO$. This is used in section~\ref{s:10} where we complete the proof of Theorem~\ref{th:1.4}. In section~\ref{s:last} we deduce that if $\J$ contains all Morse singularity types, then the loop space $\Omega^{\infty}_{\J}\MO$ is weakly homotopy equivalent to the space $\Omega B\mathcal{C}_{\J}$. 

\begin{example} The set $\Omega^{\infty-1}_{\emptyset}\MO$ is defined to be the set of maps 
\[
\R^1\times S^{\infty+d-1}\to \MO
\] 
transversal to $\BO$. Though $\Omega^{\infty-1}_{\emptyset}\MO$ is weakly homotopy equivalent to $B\mathcal{C}_{\emptyset}$, the space $\Omega^{\infty}_{\emptyset}\MO$ is not weakly homotopy equivalent to $\Omega B\mathcal{C}_{\emptyset}$.
\end{example}

Every point in $\Omega^{\infty}_{\J}\MO$ determines a fiber $F\subset S^{\infty+d}$; namely a point $f: S^{\infty+d}\to \MO$ determines the fiber $F=f^{-1}(\BO)$. In section~\ref{s:12} we show that every point in $\Omega^{\infty}_{\J}\MO$ determines not only the set $F\subset S^{\infty}$ but also an additional structure on $F$, which we use in section~\ref{s:13} to prove Theorem~\ref{th:main}.

It is interesting that in Theorem~\ref{th:1.4} we do not require that $\J$ contains all Morse singularity types, while this requirement is essential for Theorem~\ref{th:10.2}, which we use in order to establish Theorem~\ref{th:main}.

\addcontentsline{toc}{part}{Sheaves of singular manifolds}

\section{Sheaves of singular manifolds} \label{sec:3}
\subsection{Singularities of smooth maps}
A smooth map $f: M\to N$ of manifolds is said to be {\it singular} at a point $x\in M$ if the rank of the differential $df$ at $x$ is strictly less than the minimum of $\dim M$ and $\dim N$. In this case $x$ is said to be a \emph{singular point} of $f$. We are only interested in the case where the {\it dimension} $d$ of $f$, i.e., the difference $\dim M-\dim N$, is non-negative. 

Two continuous maps $f,g\colon X\to Y$ of topological spaces define the same \emph{map germ} at $x\in X$ if there is a neighborhood $U\subset X$ of $x$ such that $f|U=g|U$. 
We use the notation 
\begin{equation}\label{eqi:3.1}
    f\colon (X, x)\longrightarrow (Y, f(x))	 
\end{equation}
for the map germ at $x$ defined by $f$. In fact the same notation (\ref{eqi:3.1}) is used even if the map $f$ is defined only on a small neighborhood of $x$ in $X$. For $i=1,2$, let $X_i$ and $Y_i$ be smooth manifolds. We say that smooth map germs $f_i: (X_i, x_i)\to (Y_i, f(x_i))$ are of the same \emph{singularity type} if there are neighborhoods $U_i\subset X_i$ of $x_i$ and $V_i\subset Y_i$ of $f_i(x_i)$, and diffeomorphism germs $\alpha, \beta$ that fit a commutative diagram of map germs
\[
\begin{CD}
(U_1, x_1) @>\alpha>> (U_2, x_2) \\
@Vf_1 VV @V f_2 VV \\
(V_1, f_1(x_1)) @>\beta>> (V_2, f_2(x_2)),
\end{CD}
\]
We also say in this case that the singular points $x_1$ of $f_1$ and $x_2$ of $f_2$ are of the same singularity types. 

\begin{convention} For any singularity type $\tau$ there is a well-defined dimension $d=d(\tau)$, which is the dimension of maps that have singularities of type $\tau$. We will assume that a dimension $d\ge 0$ is chosen and fixed throughout the paper and that each set $\J$ of singularity types that we consider consists of types $\tau$ with $d(\tau)=d$. 
\end{convention}

An \emph{unfolding} of a smooth map germ $f: X\to Y$ at $x$ is a cartesian diagram of smooth map germs 
\[
\begin{CD}
(X', x') @>F >> (Y', F(x')) \\
@A iAA @AjAA \\
(X, x) @> f>> (Y, f(x)) 
\end{CD}
\]
where $X'$ and $Y'$ are smooth manifolds, $x'$ is a point in $X'$, and $i$ and $j$ are immersion germs such that $j$ is transverse to $F$. Here \emph{cartesian} means that $(f,i)$ is a diffeomorphism germ at $x$ of $X$ onto $\{\,(u,v)\in X'\times Y\,|\,F(u)=j(v) \,\}$. An unfolding is \emph{trivial} it 
there are map germs $r: X'\to X$ at $x'$ and $s: Y'\to Y$ at $F(x')$ such that $r\circ i=\id_X$, $s\circ j=\id_Y$ and $f\circ r=s\circ F$. A map germ $f$ is \emph{stable} if each of its unfoldings is trivial. An unfolding is \emph{stable} if it is stable as a map germ.

\begin{remark} It is known that the subset of map germs with no stable unfoldings is of infinite codimension in the space of all map germs. In other words every `generic' map germ has a stable unfolding. Thus we may assume that all map germs under consideration have stable unfoldings (e.g., see \cite{GWPL}). It is also known that any two stable unfoldings of the same dimension of the same map germ are isomorphic; and any stable map germ is given by a germ of a polynomial with respect to suitable coordinates (e.g., see \cite{GWPL}).  
\end{remark}

For a set $\J$ of singularity types, a smooth map $f$ is said to be a {\it $\J$-map} if each of its singular points is of type in the set $\J$. Two map germs are \emph{of the same stable singularity type} or \emph{stably equivalent}  if they have unfoldings of the same singularity type. A set $\J$ is \emph{stable} if $\tau\in \J$ implies that all singularity types stably equivalent to $\tau$ are in $\J$ as well. A set $\J$ of singularity types is said to be \emph{open} if every map $C^{\infty}$-close to a $\J$-map is a $\J$-map itself.

\subsection{The Whitney topology on the space of smooth maps}

In this subsection we recall the definition of the strong topology on the space of continuous maps $C^0(X, Y)$ of topological spaces $X$ and $Y$, and the definition of the Whitney topology on the space $C^{\infty}(M, N)$ of smooth manifolds $M$ and $N$. 

Given a continuous map $f\colon X\to Y$ of topological spaces, the base of the \emph{strong topology} on the space $C^{0}(X, Y)$ about $f$ is formed by the sets 
\[
   \mathcal{N}_{W}\colon= \{\, g\in C^0(X, Y)\, |\, \Gamma_g\subset W\,\}
\]
where $\Gamma_g\subset X\times Y$ is the graph of $g$ and $W\subset X\times Y$ is an arbitrary open set containing the graph of $f$.  

The Whitney topology on the space of smooth maps is defined by means of $k$-jets. We say that two map germs $f, g: M\to N$ at $x\in M$ represent the same \emph{$k$-jet} if $f(x)=g(x)$ and all derivatives of $f$ and $g$ at $x$ of order $\le k$ are the same with respect to some (and hence any) pair of coordinate neighborhoods about $x$ in $M$ and about $f(x)=g(x)$ in $N$. The set of all $k$-jets forms a topological space $J^k(M, N)$ called the \emph{space of $k$-jets} of maps of $M$ into $N$. For each smooth map $f: M\to N$, there is a so called \emph{$k$-jet extension} $J^k_f: M\to J^k(M, N)$ that takes a point $x\in M$ to the $k$-jet of $f$ at $x$. In particular, there is a map 
\[
     J^k: C^\infty(M, N)\longrightarrow C^0(M, J^k(M,N))
\]  
The \emph{Whitney topology} on the space $C^{\infty}(M,N)$ of smooth maps from $M$ to $N$ is defined to be the weakest topology for which $J^k$ is continuous with respect to the strong topology on $C^0(M, J^k(M,N))$ for each $k\ge 0$.

\subsection{Sheaves}

Let us recall that a smooth map $f\colon M\to N$ is a \emph{submersion} if $\dim M\ge \dim N$ and $\mathop{\mathrm{rank}}df=\dim N$ at each point of $M$. Clearly the composition of two submersions is a submersion. 

\begin{convention}
In the paper we will consider a number of categories. Often in the definition of a category we will define only morphisms; the definition of objects can be deduced from the fact that objects in a category are in bijective correspondence with identity morphisms.  
\end{convention}

\begin{definition}\label{d:3.1} A contravariant functor $\mathcal{F}$ from the category $\mathcal{E}$ of submersions of smooth manifolds without boundary to a category $\mathcal{C}$ is a {\it sheaf on $\mathcal{E}$} if for any open covering $\{U_i\}$ of any manifold $X\in \mathcal{E}$, and sections $s_i\in \mathcal{F}(U_i)$ over each $U_i$, with $s_i=s_j$ over $U_i\cap U_j$ for all $i,j$, there is a unique section $s\in \mathcal{F}(X)$ such that $s=s_i$ over $U_i$ for all $i$. 
\end{definition}

\begin{remark}\label{rem:3.1} Our Definition~\ref{d:3.1} of a sheaf is different from that in the papers \cite{MW} of Madsen and Weiss and \cite{GMTW} of Galatius, Tillmann, Madsen and Weiss; our sheaf is defined on the category of submersions, not on the category $\mathcal{X}$ of smooth maps of manifolds without boundary. The latter sheaves are defined by replacing $\mathcal{E}$ with $\mathcal{X}$ in the Definition~\ref{d:3.1} and will be called \emph{sheaves on $\mathcal{X}$}. 
\end{remark}


Given a sheaf $\mathcal{F}$ on $\mathcal{E}$ or $\mathcal{X}$, we say that two elements $s_0$ and $s_1$ in $\mathcal{F}(X)$ are \emph{concordant} if there exists an element $s\in \mathcal{F}(X\times \R)$ such that $s$ agrees with the pullback $\mathrm{pr}^*(s_0)$ on a neighborhood of $X\times (-\infty, 0]$ and with $\mathrm{pr}^*(s_1)$ on a neighborhood of $X\times [1, \infty)$, where $\mathrm{pr}$ is the projection of $X\times \R\to X$ along the second factor. 

The set of concordant classes in a set valued sheaf $\mathcal{F}(X)$ is denoted by $\mathcal{F}[X]$. For a set valued sheaf $\mathcal{F}$ on $\mathcal{X}$, there is a canonically constructed classifying space $|\mathcal{F}|$ such that the set $\mathcal{F}[X]$ is in bijective correspondence with homotopy classes of maps $X\to |\mathcal{F}|$ for every manifold $X$ without boundary (see \cite[Proposition A.1.1]{MW}). 

Each map $\mathcal{F}_1\to \mathcal{F}_2$ of sheaves on $\mathcal{X}$ induces a map $|\mathcal{F}_1|\to |\mathcal{F}_2|$ of classifying spaces. If the induced map is a weak homotopy equivalence, then the map $\mathcal{F}_1\to \mathcal{F}_2$ is said to be a \emph{weak equivalence}.

\subsection{Sheaves of singular manifolds}\label{ss:3.4}

Let $\J$ be an arbitrary set of singularity types of map germs of dimension $d$. We say that a submanifold $W$ of 
$U\times \R\times \R^{n+d-1}$ with projections $\pi, f$ and $j$ of $W$ onto the respective factors is a {\it $\J$-submanifold} if  
\begin{itemize}
\item $\pi: W\to U$ is a $\J$-map, and
\item $(\pi, f): W\to U\times \R$ is proper.  
\end{itemize}

In view of the inclusion $\R^{n+d-1}\subset \R^{n+d}$, each $\J$-submanifold of $U\times \R\times \R^{n+d-1}$ is also a $\J$-submanifold of $U\times \R\times \R^{n+d}$. In particular there is an inclusion $D_{\J}(U; n)\to D_{\J}(U; n+1)$ of sets. For a positive integer $n$ and a manifold $U$, let $D_{\J}(U; n)$ denote the set of $\J$-submanifolds of $U\times \R\times \R^{n+d-1}$. Then the set valued functor $D_{\J}$ on $\mathcal{E}$ given by 
\[
     U\mapsto \mathop{\mathrm{colim}}_{n\to \infty} D_{\J}(U, n)
\]
is a set valued sheaf on $\mathcal{E}$, which, as we will shortly see, can be described by an infinite loop space $\Omega^{\infty-1}\mathbf{B}_{\J}$ (see Theorem~\ref{th:i1}). 

\begin{remark}\label{r:3.5} We note that each morphism $f$ in $\mathcal{E}$ determines a map $D_{\J}(f)$ of sets so that $D_{\J}$ is a contravariant functor on $\mathcal{E}$. However, if $\J$ is not empty, then a smooth map $f\colon M\to N$ of manifolds does not define a map of sets $D_{\J}(N)\to D_{\J}(M)$. This is why we are forced to consider not the category of smooth maps of manifolds without boundary as in \cite{MW} and \cite{GMTW}, but the category $\mathcal{E}$ with a smaller set of morphisms. In particular, our definition of $D_{\emptyset}$ is different from that in \cite{MW} and \cite{GMTW}. 
\end{remark}

\subsection{Infinite loop space $\Omega^{\infty}\mathbf{B}_{\J}$}\label{s:1.3}

Each set $\J$ of singularity types of map germs of dimension $d$ gives rise to a spectrum $\mathbf{B}_{\J}$ defined as follows. Let $p: E_t\to \BO_t$ be the universal vector bundle of dimension $t$ over the space $\BO_t$ of vector subspaces of $\R^{\infty}$ of dimension $t$. A point $J_f$ in a space $S_t$ is represented by a map germ $f$ at $0$ of a map of $\R^{t+d}$ into $E_t$ such that 
\begin{itemize}
\item the image of $f$ belongs to a single fiber $E_t|b$ of $p$ over some point $b\in \BO_t$,
\item $f(0)$ is the zero in the vector space $E_t|b$, and
\item the map germ $f\colon (\R^{t+d}, 0)\to (E_t|b, 0)$ is smooth, and it is of type $\J$.  
\end{itemize}
Two map germs 
\[
f, g\colon (\R^{t+d}, 0)\to (E_t|b, 0)
\]
represent the same point $J_f=J_g$ in $S_t$ if derivatives of all orders of $f$ and $g$ at $0$ are the same. The space $S_t$ is endowed with an obvious topology so that the map $\theta: S_t\to \BO_t$ that takes $J_f$ onto $b$ has a structure of a fiber bundle. The $(t+d)$-th term of the spectrum $\mathbf{B}_{\J}$ is
defined to be the Thom space $T\theta^*E_t$ of the bundle $\theta^*E_t$
over $S_t$. The desired spectrum $\mathbf{B}_{\J}$ is defined to be the
Thom spectrum with terms $T\theta^*E_t$. Its infinite loop space, i.e., the colimit of spaces $\Omega^{t+d}T\theta^*E_t$, is denoted by $\Omega^{\infty}\mathbf{B}_{\J}$.

\begin{example} A map $f:M\to N$ is said to be a \emph{submersion} if it is a map of non-negative dimension $d$  and $\mathop{\mathrm{rank}}df=\dim N$ at each point of $M$. In this case the space $S_t$ is homotopy equivalent to the fiber bundle over $\BO_t$ with fiber over $b\in \BO_t$ given by the space of surjective homomorphisms $\R^{t+d}\to E_t|b$ of vector spaces. There is a fibration map $\pi$ from $S_t$ to the Grassmannian manifold $G_d(\R^{t+d})$ of subspaces of $\R^{t+d}$ of dimension $d$. It takes a point $J_f$ in $S_t$ onto the point in the Grassmannian manifold corresponding to the kernel of $df$ at $0\in \R^{t+d}$. Let $L\subset \R^{t+d}$ be a subspace of dimension $d$, it corresponds to a point in $G_d(\R^{t+d})$. Let $L^{\perp}\subset \R^{t+d}$ denote the subspace orthogonal to $L$. Then the fiber of $\pi$ over $L$ consists of all isomorphisms from $L^{\perp}$ to the fibers of $p$; hence, the fibers of $\pi$ are contractible.    
Clearly $\mathop{\mathrm{colim}} S_t=\BO_d$, and therefore in this case $\Omega^{\infty}\mathbf{B}_{\J}$ is weakly homotopy equivalent to the space $\Omega^{\infty}B_{\emptyset}$ defined in \cite{GMTW}.  
\end{example}

\begin{remark}\label{symmetry}
Let $f\colon (M,x)\to (N,f(x))$ be a stable map germ representing a singularity type $\tau$ such that the dimension of $M$ is minimal among dimensions of source manifolds of stable map germs representing $\tau$. Then a \emph{diffeomorphism} of $\tau$ is a pair $(\alpha, \beta)$ of diffeomorphism germs $\alpha$ of $(M, x)$ and $\beta$ of $(N, f(x))$ such that $f=\beta\circ f\circ \alpha^{-1}$. The group $\Diff \tau$ is defined to be the group of diffeomorphisms of $\tau$. We note that the singularity type $\tau$ of a submersion germ is represented by a stable map germ $f$ onto a point. Clearly, in this case the group $\Diff\tau$ reduces to $\mathop{\mathrm{O}}_d$. 

By Kazarian-Szucs theorem~\cite{Kaz}, \cite{Sz} the space $\mathop{\mathrm{colim}} S_t$ breaks into the union of classifying spaces $\BDiff\tau$ of diffeomorphism groups of singularity types $\tau\in \J$ and $\BO_d$. 
\end{remark}

\section{The Galatius-Madsen-Tillmann-Weiss theorem}\label{se:4}

Let us recall \cite{GMTW} the definition of the {\it cobordism category} $\mathcal{C}_d$ of embedded manifolds of dimension $d$. An object of $\mathcal{C}_d$ is a closed smooth manifold $M$ of dimension $d-1$ embedded into a horizontal infinite dimensional plane
\[
\{a\}\times \R^{\infty+d-1}\subset \R\times \R^{\infty+d-1}= \R\times\mathop{\mathrm{colim}}_{n\to \infty}\, \R^{n+d-1}. 
\] 
The set of non-trivial morphisms in $\mathcal{C}_d$ between objects  
\[
M_1\subset \{a_1\}\times \R^{\infty+d-1} \qquad \mathrm{and}\qquad   
M_2\subset \{a_2\}\times \R^{\infty+d-1},
\]
is non-empty only if $a_1< a_2$, in which case it consists of compact smooth manifolds $W$ of dimension $d$ embedded into $[a_1, a_2]\times \R^{\infty+d-1}$ such that the boundary $\partial W$ is given by $M_1\sqcup M_2$, and for some $\varepsilon>0$,
\[
W\cap ([a_0, a_0+\varepsilon)\times \R^{\infty+d-1})=[a_0, a_0+\varepsilon)\times M_0, 
\]
\[
W\cap ((a_1-\varepsilon, a_1]\times \R^{\infty+d-1})=(a_1-\varepsilon, a_1]\times M_1.  
\]
There are natural topologies on the sets of objects and morphisms of $\mathcal{C}_d$ so that the space of objects in $\mathcal{C}_d$ is weakly homotopy equivalent to the disjoint union $\sqcup \BDiff M$ of classifying spaces of diffeomorphism groups of manifolds of dimension $d-1$, while the space of non-trivial morphisms is weakly homotopy equivalent to the disjoint union 
\[
\sqcup \BDiff (W, M_1, M_2)
\] 
of classifying spaces of groups of diffeomorphisms of $W$ that restrict to diffeomorphisms of open neighborhoods of the incoming boundary $M_1$ and outcoming boundary $M_2$.  

The category $\mathcal{C}_d$ is a \emph{topological category}, i.e., a category in which the spaces of objects and morphisms are topological spaces and the structure maps of the category $\mathcal{C}_d$ (source, target, identity and composition) are continuous.

The main theorem in \cite{GMTW} relates the classifying space $B\mathcal{C}_d$ of the cobordism category and the infinite loop space $\Omega^{\infty}\mathbf{B}_{\emptyset}$ where $\emptyset$ is the empty set of singularities of maps of dimension $d$. 

\begin{theorem}[Galatius, Madsen, Tillmann, Weiss, \cite{GMTW}]\label{th:GMTW} There is a weak homotopy equivalence $B\mathcal{C}_d\to \Omega^{\infty-1}\mathbf{B}_{\emptyset}$.
\end{theorem}

In order to prove Theorem~\ref{th:GMTW}, the authors of \cite{GMTW} related  $D_{\emptyset}$ and $\mathcal{C}_d$ and proved the following theorem, which is interesting in its own right.  

\begin{theorem}[Galatius, Madsen, Tillmann, Weiss, \cite{GMTW}]\label{th:GMTW2} For any manifold $X$ there is a bijection 
\[
    D_{\emptyset}[X] \longrightarrow [X, \Omega^{\infty-1}\mathbf{B}_{\emptyset}]  
\]
between the set of concordance classes of elements in $D_{\emptyset}(X)$ and the set of homotopy classes of maps of $X$ into $\Omega^{\infty-1}\mathbf{B}_{\emptyset}$. 
\end{theorem}

In the next section we prove Theorem~\ref{th:i1} which is a generalization of Theorem~\ref{th:GMTW2}.

\section{Singular version of Theorem~\ref{th:GMTW2}}\label{s:5}

The infinite loop space $\Omega^{\infty}\mathbf{B_{\J}}$ determines a cohomology theory $H^*(-; \Omega^{\infty}\mathbf{B_{\J}})$; the $s$-th cohomology group of a CW-complex $X$ is given by 
\[
     H^s(X; \Omega^{\infty}\mathbf{B_{\J}})=[X, \Omega^{\infty-s}\mathbf{B_{\J}}].
\]
Furthermore the space $\Omega^{\infty}\mathbf{B_{\J}}$ is defined so that the cohomology group $H^{0}(X; \mathbf{B_{\J}})$ models cobordism groups of maps with $\J$-singularities (see section~\ref{s:last}). 

\begin{remark}
In general, the group $H^{0}(X; \Omega^{\infty}\mathbf{B}_{\J})$ may not be isomorphic to the cobordism group of $\J$-maps into $X$; the assumptions in the statement of Theorem~\ref{th:0.2} can not be omitted. However, in this section we show (see Theorem~\ref{th:i1}) that the group $H^{1}(X; \Omega^{\infty}\mathbf{B}_{\J})$ is precisely what one expects it to be even if $\J$ contains no Morse singularity types. 
\end{remark}

As has been mentioned, in the case of non-singular maps of dimension $d$, i.e., in the case where $\J=\emptyset$, the space $\Omega^{\infty}\mathbf{B}_{\J}$ naturally appears in the paper \cite{GMTW}; according to Theorem~\ref{th:GMTW2}, the homotopy classes of maps $[N, \Omega^{\infty-1}\mathbf{B}_{\emptyset}]$ of a closed manifold $N$ are in bijective correspondence with the concordance classes $D_{\emptyset}[N]$ of non-singular maps to $N$. In this section we extend Theorem~\ref{th:GMTW2} to the case of an open stable set $\J$ of singularity types; this is one of the steps in the proof of Theorem~\ref{th:1.4}.  

\begin{theorem}\label{th:i1} Let $\J$ be an open stable set of singularity types. Then, for each closed manifold $N$, the homotopy classes $[N, \Omega^{\infty-1}\mathbf{B}_{\J}]$ of maps of $N$ into $\Omega^{\infty-1}\mathbf{B}_{\J}$ are in bijective correspondence with the concordance classes $D_{\J}[N]$ of $\J$-maps into $N$. 
\end{theorem}

Theorem~\ref{th:i1} is proved in subsection~\ref{ss:proof}. 

\subsection{Spaces of Taylor polynomials}

Let $k$ be a positive integer or $`\infty'$. Given a smooth map $f\colon \R^m\to \R^n$ with $f(0)=0$ the \emph{$k$-jet extension} of $f$ at $0$ is defined to be the sequence 
\[
   j^kf(0)=(f'(0), f''(0), \dots, f^{(k)}(0))
\]
of all derivatives of $f$ at $0$ of order $\le k$. The space  $J^k_0(\R^m, \R^n)$ of all $k$-jet extensions of smooth maps $f\colon \R^m\to \R^n$ with $f(0)=0$ has a structure of a smooth manifold; in fact, it can be canonically identified with $\R^N$ for some positive integer $N$. For example, 
\[
    J^1_0(\R^m, \R^n)= \mathop{\mathrm{Hom}}(\R^m, \R^n)\simeq \R^{mn}
\]
where $\mathop{\mathrm{Hom}}(\R^m, \R^n)$ is the space of linear homomorphisms from $\R^m$ to $\R^n$. 

Let us now introduce a parametric version of $J^k_0(\R^m, \R^n)$. For a vector bundle $\gamma$ over a space $X$ we will denote the fiber of $\gamma$ over $x\in X$ by $\gamma_x$. Given two vector bundles $\xi$ and $\eta$ over a smooth manifold $M$ of dimensions $m$ and $n$ respectively, we define a fiber bundle 
\[
   J^k(\xi, \eta) \longrightarrow M
\]
with fiber isomorphic to $J^k_0(\R^m, \R^n)$. A point of $J^k(\xi, \eta)$ in the fiber over $x\in M$ is an equivalence class of map germs $f\colon (\xi_x, 0)\to (\eta_x, 0)$. Here two map germs $f$ and $g$ are equivalent if all derivatives of $f$ and $g$ at $0$ of order $\le k$ are the same. In order to justify the equivalence relation, let us recall that if the derivatives of order $\le k$ of two map germs are the same at some point with respect to one choice of local coordinates, then the same is true for any other choice of local coordinates. The equivalence class of a map germ $f$ is called the \emph{$k$-jet extension} of $f$; it is also denoted by $j^kf$. The space $J^k(\xi, \eta)$ is called a \emph{$k$-jet space}. 

For a manifold $X$, the total space of the tangent bundle of $X$ is denoted by $TX$. The tangent space of $X$ at a point $x$ is denoted by $T_xX$. 

Given a continuous map $f: M\to N$, we will often be interested in the $k$-jet space $J^k(TM, f^*TN)$. To simplify notation, we will denote this space by $J^k(TM, TN)$. In fact we will often write $TN$ for $f^*TN$ if the map $f$ is understood.

\subsection{Proof of Theorem~\ref{th:i1}}\label{ss:proof}
We will tacitly assume that all manifolds are equipped with a Riemannian metric; the choice of Riemannian metrics on manifolds under consideration is not important. 

\begin{lemma}\label{l:4.1} There is a well defined map $\rho\colon D_{\J}[X] \longrightarrow [X, \Omega^{\infty-1}\mathbf{B}_{\J}]$ of sets. 
\end{lemma}
\begin{proof} Let $W$ be a representative of $D_{\J}[X]$, i.e., $W$ is a submanifold of 
\[
X\times \R\times \R^{n+d-1}=X\times \R^{n+d-1}\times \R
\] 
with projections $\pi$, $f$ and $j$ of $W$ onto the respective factors of $X\times \R\times \R^{n+d-1}$. Then $\pi: W\to X$ is a $\J$-map, and, in particular, if $W$ is of dimension $m$, then the dimension of $X$ is $m-d$. Without loss of generality we may assume that $n\gg \dim W$ and $f\times j$ is an embedding of $W$. We will need several bundles over $W$. Let $\gamma$ be the normal bundle of $W$ in $\R\times \R^{n+d-1}$; it is of dimension $n+d-m$. Let $\nu$ be the normal bundle of dimension $n$ of $W$ in $X\times \R\times \R^{n+d-1}$. We note that there is a canonical isomorphism 
\[
     \nu= TX\oplus \gamma
\]  
of vector bundles over $W$. Let $\varepsilon$ be the trivial vector bundle of dimension $1$ over $W$. Given a point $x$ in $W$, we may use the Riemannian structure on $W$ to canonically identify a neighborhood of $x$ in $W$ with a neighborhood of $0$ in $T_xW$. Similarly we can identify a neighborhood of $\pi(x)$ in $X$ with a neighborhood of $0$ in $T_{\pi(x)}X$. 
Then $\pi$ gives rise to a section $s$ of the $k$-jet bundle
\[
   J^k(TW\oplus \gamma, TX\oplus \gamma)=J^k(\varepsilon^{n+d}, \nu) \longrightarrow W; 
\]
this section takes a point $x\in W$ to the $k$-jet extension of the map germ
\[
\pi\times \id|T_xW\times \gamma_x\colon T_xW\times \gamma_x\longrightarrow T_{\pi(x)}X\times \gamma_x
\] 
at $x\times 0$ where $\id$ is the identity map of the fiber $\gamma$ over $x$. Since $\pi$ is a $\J$-map and $\J$ is a stable set of singularity types, the image of the section $s$ is in the subset of $k$-jets of germs of $\J$-maps. Consequently, the section $s$ determines a lift of the map $W\stackrel{\nu}\longrightarrow \BO_{n}$ to a map $W\to S_{n}(\J)$ with respect to the projection $\theta: S_n(\J)\to \BO_n$. In other words there is a commutative diagram 

\[
\xymatrix{
  &                              &  S_n(\J) \ar[d]^{\theta} & \\
  &           W   \ar[ur] \ar[r]^{\nu} & \BO_n &
}
\]

Now the Pontrjagin-Thom construction (e.g., see \cite{St}) yields a map
\[
    X_+\wedge S^{n+d-1} \longrightarrow T\theta^*E_{n}
\]
whose adjoint represents a desired map
\begin{equation}\label{eq:4.5}
    X \longrightarrow  \Omega^{\infty-1}\mathbf{B}_{\J}.
\end{equation}
The homotopy class of the map $(\ref{eq:4.5})$ does not depend on the choice of the representative of the concordance class of $W$. Hence there is a well defined map
\[
    \rho\colon D_{\J}[X] \longrightarrow [X, \Omega^{\infty-1}\mathbf{B}_{\J}]. 
\]
\end{proof}

\begin{lemma}\label{l:4.2} There is a well defined map
    $\sigma: [X, \Omega^{\infty-1}\mathbf{B}_{\J}] \longrightarrow D_{\J}[X].$
\end{lemma}

\begin{proof}
An element of $[X, \Omega^{\infty-1}\mathbf{B}_{\J}]$ is represented by a map
\begin{equation}\label{eq:4.1}
    X\longrightarrow \Omega^{n+d-1}T\theta^*E_{n},
\end{equation}
which, by the construction reverse to the one described in the proof of Lemma~\ref{l:4.1}, gives rise to a map $M\to X$ of a manifold of dimension $m-1$ and a section $s$ of the $k$-jet bundle 
\begin{equation}\label{eqi:1}
J^k(\varepsilon^{n+d}, \nu_M)\longrightarrow M,
\end{equation}
where $\nu_M$ is the normal bundle of $M$ in $X\times \R^{n+d-1}$. Furthermore, the image of the section $s$ is in the subset of $k$-jets of map germs of $\J$-maps.  As above we obtain a section $s'$ of the $k$-jet bundle
\begin{equation}\label{eq:3.1}
J^k(TM\oplus \varepsilon^{n+d}, TX\oplus \varepsilon^{n+d-1})\to M;
\end{equation} 
for $x\in M$ if $s(x)$ is the $k$-jet represented by a map germ $f$, then the $k$-jet $s'(x)$ in (\ref{eq:3.1}) is represented by a map germ $\id_M\times f$, where $\id_M$ is the identity map of $M$.  

We note that for an arbitrarily chosen closed $(n+d-1)$-dimensional parallelized manifold $L$ the tangent bundle of $M\times \R\times L$ is canonically isomorphic to $TM\oplus \varepsilon^{n+d}$, while the tangent bundle of $X\times L$ is canonically isomorphic to $TX\oplus \varepsilon^{n+d-1}$. Thus, by the Gromov h-principle~\cite{Gr} (see also \cite{EM}) for differential relations over open manifolds, the section of (\ref{eq:3.1}) leads to a $\J$-map
\begin{equation}\label{eq:4.3}
    \alpha: M\times \R\times L\longrightarrow X\times L
\end{equation}
where $L$ is an arbitrarily chosen closed $(n+d-1)$-dimensional parallelizable manifold. 

Let $\beta$ be the projection of $M\times \R\times L$ onto the second factor. 
Let $pt\in L$ be  a regular value of the composition of the proper map 
\[
    (\alpha, \beta): M\times \R\times L \longrightarrow X\times L\times \R
\]
and the projection of $X\times L\times \R$ onto the second factor. Then for $W=(\alpha, \beta)^{-1}(X\times pt\times \R)$, the projection 
\[
    (\alpha, \beta)|W: W\longrightarrow  X\times \R 
\]
is a proper map, and, since $\J$ is a stable set of singularity types, the composition of $(\alpha, \beta)|W$ and the projection of $X\times \R$ onto $X$ is a $\J$-map.
We lift $(\alpha, \beta)|W$ to a map
\[
     W \longrightarrow X\times \R\times \R^{n+d-1}
\]
which represents an element in $D_{\J}[X]$. The obtained element in $D_{\J}[X]$ does not depend on the choice of a representative of the homotopy class of (\ref{eq:4.1}). Hence there is a well-defined map
\[
    \sigma: [X, \Omega^{\infty-1}\mathbf{B}_{\J}] \longrightarrow D_{\J}[X].
\]
\end{proof}

In order to complete the proof of Theorem~\ref{th:i1} it is necessary to show that the compositions $\sigma\circ \rho$ and $\rho\circ\sigma$ are identity homomorphisms. Proofs of both statements are left to the reader. The former is simple, while the latter is a bit technical; however, in view of the constructions in Lemmas~\ref{l:4.1} and \ref{l:4.2}, both are virtually the same as the corresponding arguments in \cite{Sa}. 

\begin{remark} In the paper~\cite{Sa} we require that the set $\J$ is $\K$-invariant. Since $\J$ is stable, this condition is satisfied (e.g., see~\cite{GWPL}). 
\end{remark}

\begin{remark} Theorem~\ref{th:i1} extends to the case of oriented $\J$-maps. In this case the sheaf $D_{\J}$ of $\J$-submanifolds is replaced by the sheaf $D^+_{\J}$ of \emph{oriented $\J$-submanifolds}. Namely, an \emph{orientation} of a $\J$-submanifold $W$ of $U\times \R\times \R^{n+d-1}$ with projections $\pi$, $f$ and $j$ of $W$ onto the respective factors is an orientation of the normal bundle of $W$ in $U\times \R\times \R^{n+d-1}$. The definition of $D^+_{\J}$ is obtained from the definition of $D_{\J}$ by replacing $\J$-submanifolds by oriented $\J$-submanfolds (see subsection~\ref{ss:3.4}). The definition of the spectrum $\mathbf{B}^+_{\J}$ is obtained from the definition of the spectrum $\mathbf{B}_{\J}$ by replacing the spaces $\BO_t$ by $\mathop{\mathrm{BSO}}_t$ and by replacing the universal vector bundles of dimension $t$ by universal oriented vector bundles of dimension $t$ (see subsection~\ref{s:1.3}). The oriented version of Theorem~\ref{th:i1} asserts that for each closed manifold $N$, the homotopy classes of maps $[N, \Omega^{\infty}\mathbf{B}^+_{\J}]$ are in bijective correspondence with the concordance classes $D^+_{\J}[N]$.   
\end{remark}

\begin{remark} Theorem~\ref{ln:2} implies that in Theorem~\ref{th:i1} the assumption that the manifold $N$ is closed is not essential; the statement of Theorem~\ref{th:i1} holds true for open manifolds as well. 
\end{remark}

\addcontentsline{toc}{part}{Singular cobordism categories}

\section{Spaces of singular manifolds} \label{s:6}

Let $\mo_n$ denote the Thom space of the universal vector bundle of dimension $n$. In particular there is an inclusion $\BO_n\subset \mo_n$ of the zero section of the universal vector bundle. The spaces $\mo_n$ with $n\ge 0$ form a spectrum $\MO$ called the Thom spectrum. 

Studying maps with prescribed multisingularities, Maxim Kazarian observed \cite{Kaz1}, \cite{Kaz2} (see also \cite{SY}) that the infinite loop space $\Omega^{\infty+d}\MO$ can be used to construct sets $\Omega^{\infty}_{\J}\MO\subset \Omega^{\infty+d}\MO$ such that a smooth proper map $f:M\to N$ of dimension $d$ is a $\J$-map if and only if the image of the associated Pontrjagin-Thom map $N\to \Omega^{\infty+d}\MO$ is in $\Omega^{\infty}_{\J}\MO$. In the rest of this section we will recall the construction of the set $\Omega^{\infty}_{\J}\MO$ with necessary modifications, and introduce a topology on the constructed sets.  

\subsection{Preliminary remarks}

To begin with let us recall that the space $\BO_n$ is defined to be the colimit of Grassmannian manifolds $G_n(\R^t)$ of subspaces of dimension $n$ in $\R^t$ with $t>n$. Similarly the total space $EO_n$ of the universal vector bundle over $\BO_n$ is the colimit of the total spaces $E_n(\R^t)$ of the canonical vector bundles of dimension $n$ over Grassmannian manifolds. The space $EO_n$ is, of course, equipped with the colimit topology, i.e.,  a subset in $EO_n$ is closed if and only if its intersection with $E_n(\R^t)$ is closed in $E_n(\R^t)$. In particular, each space $E_n(\R^t)$ is closed in $EO_n$. 

The image of every map of a compact set $K$ to $EO_n$ belongs to some finite dimensional manifold $E_n(\R^t)$. This can be proved by a standard argument, which we recall for the convenience of the reader. Indeed, for each $t$ we may choose a point $x_t\in K$ that maps to $E_n(\R^t)\setminus E_{n-1}(\R^t)$. Since the intersection of the image of any subset $X'$ of $X=\{x_t\}$ with every set $E_n(\R^t)$ is finite, we conclude that the image of $X'$ is closed in $EO_n$ in the colimit topology. Hence $X'$ is closed in $K$, and therefore every subset in $X$ is closed. Consequently $X$ is a discreet subset of a compact set. Thus $X$ is finite, which is equivalent to the desired statement.

\subsection{Definition of the set $\Omega^{n-k}_{\J}\mo_n$}\label{si:1}

In this subsection we define a set $\Omega^{n-k}_{\J}\mo_n$  for each non-negative integer $k$ and $n\ge k$. We will turn it into a topological space in subsection~\ref{s:6.3}.

\begin{definition}\label{d:1} We say that a map $f: M\to \mo_n$ of a smooth manifold $M$ is \emph{smooth near $f^{-1}(\BO_n)$} if there is an open neighborhood $V$ of $\BO_n$ in $\mo_n$ such that for each point $x$ in $f^{-1}(V)$ there is an open neighborhood $U$ of $x$ in $M$ over which $f|U$ is a smooth map to $E_n(\R^t)\subset \mo_n$.
\end{definition}

Let $\J$ be an open stable set of singularity types of map germs of dimension $d$.
Let $k$ be a non-negative integer and $n$ be an integer with $n>k$. 

\begin{definition} \label{d:2}
The space $\Omega^{n-k}_{\J}\mo_n$ is defined to be the set of maps
\begin{equation}\label{eq:6.1}
     f: \R^k\times S^{n+d-k}\longrightarrow {\mo}_{n},  
\end{equation}
subject to the following conditions 

\begin{itemize}
\item $f$ is smooth near $f^{-1}(\BO_n)$,
\item for a distinguished point $*$ in $S^{n+d-k}$, the map $f$ takes $\R^k\times \{*\}$ onto the distinguished point of $\mo_n$, and 
\item for each point $x$ in $f^{-1}(\BO_n)$, the composition
\[
     (\R^k\times S^{n+d-k}, x) \longrightarrow ({\mo}_n, f(x)) \longrightarrow (\R^n, 0)
\]
of the map germ $f$ at $x$ and a projection germ on a fiber of the universal vector bundle $EO_n\to \BO_n$ at $f(x)$ is a map germ of type in $\J$. 
\end{itemize}
\end{definition}

For a smooth map $f: M\to N$ of dimension $d$, the Pontrjagin-Thom construction results in a map 
\begin{equation}\label{eq:6.4}
        N\longrightarrow \Omega^{n+d}{\mo}_n,
\end{equation}
provided that $n\gg \dim M$. Kazarian~\cite{Kaz1}, \cite{Kaz2} made an essential observation that the map $f$ is a $\J$-map if and only if the image of $N$ under the map (\ref{eq:6.4}) is in the set $\Omega^n_{\J}{\mo}_n$. 

\begin{example} Suppose that the set of singularity types $\J$ is empty. Then the set $\Omega^n_{\J}{\mo}_n$ consists of maps $f: S^{n+d}\to \mo_n$ transversal to $\BO_n$. In this case the set $\Omega^n_{\J}\mo_n$ is an approximation of the disjoint union $\sqcup \BDiff M$ of classifying spaces of diffeomorphism groups of manifolds of dimension $d$ considered as a set; the greater the $n$, the better the approximation. The approximating map of $\Omega^n_{\J}\mo_n$ to $\sqcup \BDiff M$ takes a point in $\Omega^n_{\J}\mo_n$ corresponding to a map $f$ onto the point in $\sqcup \BDiff M$ corresponding to the manifold $f^{-1}(\BO_n)$ embedded into $S^{n+d}$.

We note, however, that if the topology on $\Omega^n_{\J}\mo_n$ is chosen  to be the one inherited from the compact-open topology on the mapping space $\Omega^{n+d}\mo_n$, then a point in the component of $\Omega^{n}_{\J}\mo_n$ corresponding to $\BDiff M$ is not disjoint from points in other components corresponding to $\BDiff M'$ because an arbitrarily small $C^0$-deformation of $f$ near $M=f^{-1}(\BO_n)$ may result in attachment of small handles to $M$. 
\end{example}

\subsection{Topology on the sets $\Omega^{n-k}_{\J}T$}\label{si:6.2}

Let us recall that the space $\mo_n$ is the colimit of Thom spaces of canonical vector bundles $E_n(\R^t)$ over Grassmannian manifolds for $t\ge 0$. In this subsection we will keep $n$ fixed and denote the Thom space of $E_n(\R^t)$ simply by $T=T_t$. We note that the construction in subsection~\ref{si:1} of the set $\Omega^{n-k}_{\J}\mo_n$ can be modified by replacing $\mo_n$ with $T$ to produce a set $\Omega^{n-k}_{\J}T$ of maps
\[     
            f: \R^k\times S^{n-k} \longrightarrow T 
\]
for which there is a special neighborhood $V$ of $G_n(\R^t)$ in $T$ that satisfies requirements listed in Definition~\ref{d:2}. Let $V'$ be another (smaller) neighborhood of $G_n(\R^t)$ in $T$ whose closure is contained in $V$, and let $U$ denote a neighborhood of the graph of $f$ in $\R^k\times S^{n-k}\times T$. 

\begin{remark} The idea behind our choice of topology is to control the jets of $f$ at points in $f^{-1}(V')$ and the values of $f$ elsewhere. 
\end{remark}

We will define a topology on $\Omega^{n-k}_{\J}T$ similarly to the Whitney topology by means of locally trivial fiber bundles 
\[
   \pi\colon J^s(\R^k\times S^{n-k}, E_n(\R^t))\longrightarrow \R^k\times S^{n-k}\times E_n(\R^t)
\]
for $s=1,2,...$, where the total space is the jet space of smooth maps of $\R^k\times S^{n-k}$ into $E_n(\R^t)$, and the projection $\pi$ takes the jet of a map germ $g$ at $x$ onto a point $x\times g(x)$.  Let $O$ be an open set in the jet space that contains $J^s_f(f^{-1}(V))$. Then, by definition, the base of open neighborhoods in our topology on $\Omega^{n-k}_{\J}T$ consists of subsets
\[
    \mathcal{N}(V', O, s)\colon= \{\, g\in \Omega^{n-k}_{\J}T\, |\, J^s_g(f^{-1}(V'))\subset O, \Gamma_g\subset U   \,\}. 
\]
In particular, in this topology a map $g$ is close to $f$ means that 
\begin{itemize} 
\item for each point $x$ in the domain of $f$ with $f(x)$ sufficiently close to $\BO_n$, the jets of $g$ and $f$ at $x$ are close; and
\item for each point $x$ in the domain of $f$ with $f(x)$ sufficiently far from $\BO_n$, the values $f(x)$ and $g(x)$ are close. 
\end{itemize}   

\subsection{Topology on the set $\Omega^{n-k}_{\J}\mo_n$}\label{s:6.3} For each $t\ge 0$, there are obvious inclusions of sets
\begin{equation}\label{eqi:2}
  \Omega^{n-k}_{\J}TE_n(\R^t) \longrightarrow \Omega^{n-k}_{\J}TE_n(\R^{t+1}).
\end{equation}
In the case where $k=0$, the colimit of inclusions (\ref{eqi:2}) of sets is precisely the set $\Omega^{n-k}_{\J}\mo_n$. Thus we may define the topology on the latter space to be the colimit topology. In the general case more care is necessary as the set $\Omega^{n-k}_{\J}\mo_n$ contains elements that are not present in the colimit of inclusions (\ref{eqi:2}).

We say that two map germs $f, g: M\to EO_n$ define the same \emph{$s$-jet} at a point $x$ of a smooth manifold $M$ if there is a neighborhood $O$ of $x$ in $M$ such that the images $f(O)$ and $g(O)$ belong to the manifold $E_n(\R^t)$ and 
\[
    f, g|O\colon O\longrightarrow E_n(\R^t)
\]
represent the same $s$-jet. The space of $s$-jets of map germs of $M$ into $EO_n$ form a topological space denoted by $J^s(M, EO_n)$. It is the total space of a fiber bundle
\[
     J^s(M, EO_n)\longrightarrow M\times EO_n.
\]
Topology on $\Omega^{n-k}_{\J}\mo_n$ is defined as in subsection~\ref{si:6.2} by replacing $E_n(\R^t)$ with $EO_n$.

\subsection{The topological space $\Omega^{\infty}_{\J}\MO$}
Finally, we turn to the definition of the topology on $\Omega^{n-k}_{\J}\mo_n$. For each $k\ge 0$ and $n>k$, there is an inclusion 
\begin{equation}\label{eq:6.2}
    \Omega^{n-k}_{\J}{\mo}_n \longrightarrow \Omega^{n-k+1}_{\J}{\mo}_{n+1}, 
\end{equation}
which is described as follows. Let $f$ be a map as in (\ref{eq:6.1}). We may regard $f$ to be a family of pointed maps $f_p$ from $S^{n+d-k}$ to $\mo_n$ parametrized by points $p\in \R^k$. Then the inclusion (\ref{eq:6.2}) takes the family of maps $f_p$ to the family of maps
\[
   S^{n+d-k+1}\longrightarrow S^1\wedge S^{n+d-k}\stackrel{\id_{S^{1}}\wedge f_p}\longrightarrow S^1\wedge {\mo}_{n}\longrightarrow {\mo}_{n+1},
\]
where $\id_{S^1}$ is the identity map of the circle $S^1$, and the last map is the structure map of the spectrum $\MO$. A straightforward verification shows that if $f$ is a point in $\Omega^{n-k}_{\J}\mo_n$, then the resulting map is a point in $\Omega^{n-k+1}_{\J}\mo_{n+1}$. The {\it set $\Omega^{\infty-k}_{\J}\MO$}, with $k\ge 0$, is defined to be the colimit of sets $\Omega^{n-k}\mo_n$. We equip it with topology defined by the colimit (\ref{eq:6.2}). The topological space $\Omega^{\infty}_{\J}\MO$ is called the {\it space of singular manifolds with $\J$-singularities}. 

\begin{remark} Spaces $\Omega^{\infty-k}_{\J}\MSO$ of oriented singular manifolds with $\J$-singularities are defined similarly by replacing in the definition of $\Omega^{\infty-k}_{\J}\MO$ the spectrum $\MO$ by the spectrum $\MSO$.  
\end{remark}

\section{Extension of $D_{\J}$}\label{s:7}

In this section we define a modification $E_{\J}$ of the set valued sheaf $D_{\J}$ on the category $\mathcal{E}$ so that for each manifold $X$, there is an identification of concordance classes $D_{\J}[X]\to E_{\J}[X]$, and, on the other hand, there is an extension of the sheaf $E_{\J}$ to the category $\mathcal{X}$ of  smooth maps of smooth manifolds without boundary.   

In fact, we define $E_{\J}$ to be the set valued sheaf on $\mathcal{X}$ that associates the set of continuous maps 
\begin{equation}\label{eq:le1}
       X \longrightarrow \Omega^{\infty-1}_{\J}\MO
\end{equation}
to each smooth manifold $X$ without boundary, where the topology on the space $\Omega^{\infty-1}_{\J}\MO$ is the one defined in section~\ref{s:6}. A smooth map $f\colon Y\to X$ gives rise to a map of sets $E_{\J}(X)\to E_{\J}(Y)$ that precomposes a map (\ref{eq:le1}) with $f$. The restriction sheaf $E_{\J}|\mathcal{E}$ can be related to $D_{\J}$ by means of the Pontrjagin-Thom construction which yields a map $D_{\J}[X]\xrightarrow{e_X} E_{\J}[X]$ of concordance classes for any manifold $X\in \mathcal{E}$. We will show that the map $e_X$ is an isomorphism. 

\begin{lemma}\label{ln:1} Suppose that $\J$ is an open stable set of singularity types. Then each element of $E_{\J}[X]$ admits a representation by a map
\[
      P: X\times \R\times S^{n+d-1} \longrightarrow   {\mo}_n 
\]
which is smooth in a neighborhood of $P^{-1}(\BO_n)$ and transversal to $\BO_n$. 
\end{lemma}
\begin{remark} The definition of a smooth map into $\mo_n$ is given in section~\ref{s:6} (see Definition~\ref{d:1}). 
\end{remark}
\begin{proof} Each element in $E_{\J}(X)$ is represented by a continuous map
\begin{equation}\label{eq:8.2}
      Q: X\times \R\times S^{n+d-1} \longrightarrow   {\mo}_n
\end{equation}
such that 
\begin{itemize}
\item for each point $x\in X$, the restriction $Q_x$ of $Q$ to the slice $\{x\}\times \R\times S^{n+d-1}$ is smooth near $Q_x^{-1}(\BO_n)$, 
\item $X\times \R\times \{*\}$ is taken onto the distinguished point in $\mo_n$, and
\item for each $y\in Q_x^{-1}(\BO_n)$ the composition of the map germ $Q_x$ at $y$ and a projection
\[
   ({\mo}_n, Q_x(y)) \longrightarrow (\R^n, 0)
\]
is a map germ of a type in $\J$. 
\end{itemize}

  There exists a $C^0$-slight perturbation $P$ of $Q$ relative to $X\times \R\times \{*\}$ such that in addition to the three properties enjoyed by $Q$, the map $P$ is smooth near $P^{-1}(\BO_n)$ and transversal to $\BO_n$. In fact we may choose $P$ so that for each $x\in X$ the restriction $P_x$ of $P$ to $\{x\}\times \R\times S^{n+d-1}$ is $C^{\infty}$-close near $Q_x^{-1}(\BO)$ to the restriction $Q_x$ of $Q$. Then, since the set $\J$ of singularity types is open, for each $x$ the map $P_x$ represents a loop in $\Omega^{\infty-1}_{\J}\mo_n$. Consequently $P$ represents a point in $E_{\J}(X)$. Furthermore, in view of the chosen topology on $\Omega^{\infty-1}_{\J}\MO$, the points represented by $P$ and $Q$ in $E_{\J}(X)$ are in the same concordance class $E_{\J}[X]$. Consequently, the map $P$ satisfies the requirements of Lemma~\ref{ln:1}.
\end{proof}

In fact a relative version of Lemma~\ref{ln:1} holds true. Namely, let $A$ be a closed subset in $X$, and $Q$ be a representative of a given element in $E_{\J}(X)$ such that $Q|A$ is smooth near $(Q|U)^{-1}(\BO_n)$ and transversal to $\BO_n$. Then the perturbation $P$ in Lemma~\ref{ln:1} can be chosen relative a compact neighborhood of $A$.

\begin{remark} If a continuous map (\ref{eq:le1}) is represented by a map (\ref{ln:1}) transversal to $\BO_n$, then it corresponds to a smooth manifold given by $Q^{-1}(\BO_n)$. For a general continuous map (\ref{eq:le1}) there is no such a nice geometric correspondence. We may restrict our attention to continuous maps (\ref{eq:le1}) represented by maps (\ref{ln:1}) transversal to $\BO_n$. However, in this case we get only a sheaf on $\mathcal{E}$, not on $\mathcal{X}$.  
\end{remark}

\begin{lemma}\label{ln:4} The space $\Omega^{\infty-1}\mathbf{B}_{\J}$ has the homotopy type of a CW-complex provided that $\J$ is an open set of singularity types of map germs of dimension $d$.
\end{lemma}
\begin{proof} To prove the lemma we follow the definition of the space $\Omega^{\infty}\mathbf{B}_{\J}$. It is the infinite loop space of a spectrum $\mathbf{B}_{\J}$. Each space of the spectrum $\mathbf{B}_{\J}$ is the Thom space of a vector bundle over a space $S_t$ which, in its turn, is the total space of a locally trivial fiber bundle $\theta=\theta_t$ over $\BO_t$. Since the set $\J$ of singularity types is open, the fiber of $\theta$ is given by an open subset of a Euclidean space. In particular, the fiber of $\theta$ has the homotopy type of a CW-complex. Thus both the fiber and the base of the locally trivial bundle $\theta$ have the homotopy type of a CW-complex. Consequently, each space $S_t$ has the homotopy type of a CW-complex and the same is true for each space $T\theta^*E_t$ in the spectrum $\mathbf{B}_{\J}$. 

The infinite loop space under consideration is defined to be the colimit of spaces $\Omega^{t-1}T\theta^*E_t$ under certain inclusions. By the Milnor theorem~\cite{Mi}, the loop space of a CW-complex has the homotopy type of a CW-complex. Thus each of the spaces $\Omega^{t-1}T\theta^*E_t$ has a homotopy type of a CW-complex. Finally, the colimit of spaces with respect to inclusions is homotopy equivalent to the telescope of inclusions. Therefore, since the spaces in the colimit have homotopy types of CW-complexes, the colimit also has a homotopy type of a CW-complex. 

This implies that the colimit $\Omega^{\infty-1}\mathbf{B}_{\J}$ has a homotopy type of a CW-complex. 
\end{proof}

\begin{theorem}\label{ln:2} Suppose that $\J$ is an open stable set of singularity types of map germs of dimension $d$. Then for each manifold $X\in \mathcal{X}$, there is an isomorphism of concordance classes $E_{\J}[X]\approx D_{\J}[X]$. In fact there is a weak homotopy equivalence 
\[
\psi: \Omega^{\infty-1} \mathbf{B}_{\J}\approx \Omega^{\infty-1}_{\J}\MO.
\] 
\end{theorem}
\begin{proof} 
By Lemma~\ref{ln:4}, the space $\Omega^{\infty-1} \mathbf{B}_{\J}$ has the homotopy type of a CW-complex. Slightly abusing notation, we will assume that $\Omega^{\infty-1}\mathbf{B}_{\J}$ is itself a CW-complex. Let 
\[
     K^0\subset K^1\subset K^2\subset \cdots \subset \Omega^{\infty-1} \mathbf{B}_{\J}.
\]
be a filtration by skeleta. 

\begin{lemma} There is an increasing sequence $\{l_i\}$ of positive integers indexed by $i\ge 0$ and inclusions $K^n\subset U_n$ of $n$-th skeleta to smooth open manifolds of dimension $l_n$ such that each $U_n$ is a closed subset of $U_{n+1}$, and there are deformations $\alpha_i: U_i\to K^i$ with  $\alpha_i|U_{i-1}=\alpha_{i-1}$ for each $i$. 
\end{lemma}
\begin{proof} We will define a number $l_n$ and an open manifold $U_n$ by induction in $n$. In fact we will define $U_n$ so that $U_n$ is an open submanifold of $\R^{l_n}$ for each $n$. For $n=0$ we choose $l_n$ to be $1$ and put $U_n=K^0\times (0,1)$. Then $K^0$ is identified with a copy $K^0\times \{0\}$ in $U_0$. Suppose that for some $n$ the number $l_{n-1}$ has been defined and $U_{n-1}$ has been constructed together with an inclusion $U_{n-1}\subset \R^{l_{n-1}}$. The $n$-th skeleton $K_n$ is obtained from $K_{n-1}$ by attaching mapping cones $C$ of maps $S^{n-1}\to K_{n-1}$. Choose an embedding of $K_n$ to $\R^{l_n}$ for $l_n\gg l_{n-1}$ extending the embedding $K_{n-1}\subset U_{n-1}\subset \R^{l_{n-1}}$ so that each mapping cone approaches $\R^{l_{n-1}}$ transversally and the intersection of $U_{n-1}$ with the interior of $C$ is empty. For sufficiently small $\delta>0$, we put 
\[   
     U'_{n-1}\colon = U_{n-1}\times (-\delta, \delta)^{l_n-l_{n-1}}\subset \R^{l_{n-1}}\times \R^{l_{n}-l_{n-1}}
\]
Then $U'_{n-1}$ is a neighborhood of $K_{n-1}$ in $\R^{l_n}$ and 
\[
U'_{n-1}\cap \R^{l_{n-1}}=U_{n-1}. 
\]
For sufficiently small $\varepsilon>0$ we attach to $U_{n-1}'$ an $\varepsilon$-neighborhood of each mapping cone $C$, one for each mapping cone $C$, so that the obtained open manifold $U_n$ contains $K_n$ and $U_{n-1}$ is closed in $U_n$. This completes the induction step. 

The existence of deformations $\alpha_i$ is obvious. 
\end{proof}

By Theorem~\ref{th:i1}, for each $i$, the composition of the retraction $\alpha_i$ and the inclusion $K^i\to \Omega^{\infty-1} \mathbf{B}_{\J}$ gives rise to an element $\gamma_i$ in $D_{\J}[U^i]$. We observe that $\gamma_i|U'_{i-1}$ is the same element as the pullback of $\gamma_{i-1}$ with respect to the projection $p_{i-1}$ of $U'_{i-1}$ onto the first factor $U_{i-1}$.

By the Pontrjagin-Thom construction each element $\gamma_i$ leads to a map $\beta_i: U_i\to \Omega^{\infty-1}_{\J}\MO$. Furthermore, since $\gamma_{i}|U'_{i-1}=p_{i-1}^*\gamma_{i-1}$, the restriction of $\beta_i$ to $U_{i-1}$ is homotopic to $\beta_{i-1}$. Consequently, we may assume that $\beta_i$ coincides with $\beta_{i-1}$ over $U_{i-1}$. Then the restrictions $\{\beta_i|K^i\}$ determine a map 
\[
\psi: \Omega^{\infty-1} \mathbf{B}_{\J}\longrightarrow \Omega^{\infty-1}_{\J}\MO.
\] 
Let us now show that the map $\psi$ induces an isomorphism $\psi_*$ of homotopy groups. 

Let $f$ be a pointed map of a sphere $S$ representing an element in the homotopy group of $\Omega^{\infty-1} \mathbf{B}_{\J}$. Since $S$ is compact, its image is in $K^n$ for some $n$. Thus the construction of $\psi$ yields a map $\psi\circ f$ that is smooth near $(\psi\circ f)^{-1}(\BO)$. Suppose that the pointed map $\psi\circ f$ extends to a map of a disc $D\supset S$. We may assume that the extension restricted to a collar neighborhood $S\times [0,\varepsilon)$ of $S$ in $D$ is given by the composition of the projection of $S\times [0,\varepsilon)$ onto the first factor and $\psi\circ f$. 
Then by Lemma~\ref{ln:1} we may perturb the extended map relative to a compact neighborhood of the boundary $\partial D=S$ so that it is smooth near the inverse image of $\BO$. Thus we get a $\J$-manifold over $D$, which by the construction in the proof of Lemma~\ref{l:4.1} defines a homotopy of $f$ to a constant map. This shows that $\psi_*$ is injective. A similar argument shows that $\psi_*$ is surjective. 
\end{proof}

\section{Singular cobordism categories $\mathcal{C}_{\J}$}\label{s:9}

In this section we give a definition of the cobordism category $\mathcal{C}_{\J}$ and introduce topologies on the spaces of objects and morphisms of $\mathcal{C}_{\J}$ so that $\mathcal{C}_{\J}$ is a topological category. 

The objects of the {\it cobordism category $\mathcal{C}_{\J}$ of singular manifolds} are similar to the objects of the topological cobordism category $\mathcal{C}_d$. More precisely, an object of $\mathcal{C}_{\J}$ is a closed submanifold
\begin{equation}\label{eq:9.1}
    M \hookrightarrow \{a\} \times S^{\infty+d-1}\times \{r\}\subset \R\times  S^{\infty+d-1}\times \R, 
\end{equation}
that has a closed tubular neighborhood of radius $r$. Thus the set $\Ob\mathcal{C}_{\J}$ is defined so that there is an obvious map of sets $\Ob\mathcal{C}_{\J}\to \mathcal{C}_d$ each fiber of which is an interval. Indeed if $M\hookrightarrow \{a\}\times S^{\infty+d-1}$ is a closed manifold representing a point in $\Ob\mathcal{C}_d$, then there is the maximal number $r'$ such that for each $r\in (0, r')$, there is a unique manifold (\ref{eq:9.1}) representing an element in the set $\Ob\mathcal{C}_{\J}$. 

To introduce a topology on the set of objects of $\mathcal{C}_{\J}$, we identify $\Ob\mathcal{C}_{\J}$ with a subspace of $\Omega^{\infty-1}_{\J}\MO\times \R$ by means of the canonical Pontrjagin-Thom construction. Namely, let 
\[
       W=\R \times M \longrightarrow \R\times S^{\infty+d-1}
\]
be the inverse image of $M$ under the obvious projection 
\[
    \R\times S^{\infty+d-1}=\R\times S^{\infty+d-1}\times\{r\} \longrightarrow \{a\}\times S^{\infty+d-1}\times  \{r\}. 
\]
Then the canonical Pontrjagin-Thom map
\[
     P: \R\times S^{\infty+d-1} \longrightarrow \MO
\]
takes $W$ to $\BO$ and the complement to the tubular neighborhood of $W$ of radius $r$ onto the distinguished point of $\MO$. It follows that $P$ represents a point $p(M)$ in $\Omega^{\infty-1}_{\J}\MO$. We claim that the map of sets  
\begin{equation}\label{eq:10.1}
      \Ob\mathcal{C}_{\J} \longrightarrow \Omega^{\infty-1}_{\J}\MO\times \R,
\end{equation}
\[
     M   \mapsto    (p(M), a) 
\]
is a bijection onto the image. Indeed, given a point $(p(M), a)$ in the image, the radius $r$ can be identified with the minimum distance in $\R\times S^{\infty+d-1}$ from $P^{-1}(\BO)$ to $P^{-1}(*)$ where $*$ is the distinguished point in the spectrum $\MO$. The corresponding element in $\Ob \mathcal{C}_j$ is 
\[
   P^{-1}(\BO)\times \{r\}\cap \{a\}\times S^{\infty+d-1}\times \{r\}.
\]
We introduce a topology on $\Ob\mathcal{C}_{\J}$ so that the map (\ref{eq:10.1}) is an embedding.

Intuitively, a non-trivial morphism $W$ in $\mathcal{C}_{\J}$ between 
\[
       M_0\hookrightarrow \{a_0\} \times S^{\infty+d-1} \times \{r_0\}
\]
and
\[
       M_1\hookrightarrow \{a_1\} \times S^{\infty+d-1}\times \{r_1\},
\]
with $a_0<a_1$, is defined to be a subspace 
\begin{equation}\label{eq:8.1}
       W \hookrightarrow [a_0, a_1] \times S^{\infty+d-1}\times \{r_0\}\times \{r_1\}
\end{equation}
such that for some $\varepsilon>0$

\begin{itemize}
\item $W\cap ([a_0, a_0+\varepsilon)\times  S^{\infty+d-1}\times \{r\}) = [a_0, a_0+\varepsilon)\times M_0\times \{r\}$,
\item $W\cap ((a_1-\varepsilon, a_1]\times S^{\infty+d-1}\times \{r\}) = (a_1-\varepsilon, a_1]\times M_1\times \{r\}$, and
\item there is a proper $\J$-map $X\to Y$ of manifolds with fiber $W$ such that $\partial X\cap W=M_1\sqcup M_2$.  
\end{itemize}

More precisely we define a non-trivial morphism between $M_1$ and $M_2$ to be a point  
\[
    q: \R\times S^{\infty+d-1} \longrightarrow \MO
\]
in
\[
\Omega^{\infty-1}_{\J}\MO\times \{r_0\}\times \{r_1\}\times \{a_0\}\times \{a_1\}= \Omega^{\infty-1}_{\J}\MO
\]   
that coincides with $p(M_0)$ in a neighborhood of $(-\infty, a_0+\varepsilon] \times S^{d-1+\infty}$ and with $p(M_1)$ in a neighborhood of $[a_1-\varepsilon, \infty) \times S^{d-1+\infty}$ for some $\varepsilon>0$. Of course, here 
\[
W=q^{-1}(\BO)\quad \cap \quad  [a_0, a_1]\times  S^{d-1+\infty} \times \{r_0\}\times \{r_1\}.
\] 
We introduce a topology on the space of morphisms $W$ so that the map
\begin{equation}\label{eq:10.3}
      \Mor\mathcal{C}_{\J} \longrightarrow 
      \Omega^{\infty-1}_{\J}\MO\times \R\times \R \times \R^{2}_+,
\end{equation}
defined on the non-trivial morphisms to be
\[
         W \mapsto (q, r_0, r_1, a_0, a_1),
\]
is an embedding. Here $\R^{2}_+$ is the open half plane of points $(a_0, a_1)$ with $a_0<a_1$.

\begin{remark} Our study of singular manifolds in sections~\ref{s:12} and \ref{s:13} shows that a non-trivial morphism is equivalent to a subspace $W$ as in (\ref{eq:8.1}) equipped with an additional structure (see Theorem~\ref{th:6.2} and the construction in the beginning of section~\ref{s:13}). 
\end{remark}

\begin{remark} Oriented singular cobordism categories $\mathcal{C}^+_{\J}$ are defined similarly to singular cobordism categories $\mathcal{C}_{\J}$. The definition of morphisms in $\mathcal{C}^+_{\J}$ is obtained from the definition of morphisms in $\mathcal{C}_{\J}$ by replacing the space $\Omega^{\infty-1}_{\J}\MO$ by the space $\Omega^{\infty-1}_{\J}\MSO$.  
\end{remark}

\section{Homotopy types of singular cobordism categories}\label{s:10}

In this section we prove Theorem~\ref{th:1.4} generalizing the main theorem in \cite{GMTW}.

Theorem~\ref{th:1.4} follows from a serious of lemmas. In this section we will introduce several sheaves and prove weak homotopy equivalences:
\[
    \Omega^{\infty-1}\mathbf{B}_{\J} \longleftarrow  |E_{\J}| \longleftarrow |\beta E^{\pitchfork}_{\J}| \simeq B|E^{\pitchfork}_{\J}| \longleftarrow B|C_{\J}|\longrightarrow B\mathcal{C}_{\J}
\]

The weak homotopy equivalence on the left hand side follows from Theorem~\ref{ln:2} and the uniqueness of the classifying space $|E_{\J}|$ up to weak homotopy equivalence.

\subsection{Category valued sheaf $E^{\pitchfork}_{\J}$}

By definition the set valued sheaf $E^{\pitchfork}_{\J}$ on $\mathcal{X}$ is a subsheaf of $E_{\J}\times C^0(-, \R)$; its value on a manifold $X$ is a pair $(f, a)$  of a continuous map 
\begin{equation}\label{eqi:10.1}
    f: X\longrightarrow \Omega^{\infty-1}_{\J}\MO  
\end{equation}
and a continuous function $a: X\to \R$ such that for each point $x\in X$, the restriction of a representative 
\[
      \R \times S^{n+d-1} \longrightarrow {\mo}_n
\]
of $f(x)$ to $\{a(x)\}\times S^{n+d-1}$ is transversal to $\BO_n$. The sheaf $E^{\pitchfork}_{\J}$ is also a category valued sheaf. Indeed, the value of $E^{\pitchfork}_{\J}$ on a manifold $X$ is a poset of pairs $(f, a)$; here 
\begin{equation}\label{eq:10.2}
    (f_0, a_0)\le (f_1, a_1)
\end{equation}
if and only if $f_0=f_1$, $a_0\le a_1$ and the set $(a_0-a_1)^{-1}(0)$ is open in $X$. The category $E^{\pitchfork}_{\J}(X)$ is the category of the described poset; it has one object for each pair $(f, a)$, and one morphism 
\[
    (f_0, a_0)\mapsto (f_1, a_1)
\]
for each inequality (\ref{eq:10.2}). 

\subsection{Set valued sheaf $\beta E^{\pitchfork}_{\J}$}
Let us recall that the {\it cocycle sheaf} (see \cite{GMTW}) of $E^{\pitchfork}_{\J}$ is defined to be the set valued sheaf $\beta E^{\pitchfork}_{\J}$ on $\mathcal{X}$ of pairs of continuous maps (\ref{eqi:10.1}) and broken functions $a$. More precisely, an element in $\beta E^{\pitchfork}(X)$ is a triple  $(\mathcal{U}, \Psi, \Phi)$ of 
\begin{itemize}
\item a locally finite open cover $\mathcal{U}=\{U_j| j\in J\}$ of $X$ indexed by a fixed uncountable set $J$, 
\item a set $\Psi\in \{c_R\}$ of elements in $E^{\pitchfork}_{\J}(U_R)$, one for each non-empty finite subset $R\subset J$, where $U_{R}=\cap_{j\in R} U_j$, and
\item a collection of morphisms $\Phi=\{\varphi_{RS}\}$ in $E^{\pitchfork}_{\J}(U_{S})$, where $\varphi_{RS}$ is a morphism from $c_S$ to $c_R|U_S$, one for each non-empty finite pair $R\subseteq S\subset J$
\end{itemize} 

such that  

\begin{itemize}
\item $\varphi_{RR}$ is the identity morphism of $c_R$, 
\item $\varphi_{RT}=(\varphi_{RS}|U_T)\circ \varphi_{ST}$ for all $R\subseteq S\subseteq T$ of finite non-empty subsets of $\J$.
\end{itemize}

It follows that for a manifold $X$ in $\mathcal{X}$ a point $(\mathcal{U}, \Psi, \Phi)$ in $E^{\pitchfork}(X)$ is a pair of a continuous map (\ref{eqi:10.1}) and a \emph{broken function} $a$ given by 
the collection of functions in $\Psi$.

\subsection{Weak equivalence $E_{\J}\leftarrow \beta E^{\pitchfork}_{\J}$ of set valued sheaves}

We will need a so-called surjectivety criteria, which we recall next. 

Let $\mathcal{F}$ be a set valued sheaf on $\mathcal{X}$. Let $A$ be a closed subset of $X\in \mathcal{X}$. A \emph{germ} at $A$ is an element of the set $\mathop{\mathrm{colim}}_{U}\mathcal{F}(X)$ where $U$ ranges over open subsets of $X$ containing $A$. Let $\mathcal{F}(X, A; s)$ be the set of elements in $\mathcal{F}(X)$ whose germs at $A$ coincide with a given germ $s$. A concordance between elements in $\mathcal{F}(X, A; s)$ is a concordance in $\mathcal{F}(X)$ that restricts to a trivial concordance of the germ $s$. The \emph{surjectivety criteria} asserts that a map $\tau:\mathcal{F}_1\to \mathcal{F}_2$ of set valued sheaves on $\mathcal{X}$ is a weak equivalence if it induces a surjective map 
\[
     \mathcal{F}_1[X, A; s]\longrightarrow \mathcal{F}_2[X, A; \tau(s)]
\]
for all $X$, $A$ and $s$.

\begin{lemma}\label{ln:3} The forgetful map $\alpha\colon \beta E^{\pitchfork}_{\J}\to E_{\J}$ is a weak equivalence of set valued sheaves. 
\end{lemma}
\begin{proof} In view of Lemma~\ref{ln:1}, the proof of Lemma~\ref{ln:3} is very similar to the proof of \cite[Proposition 4.2]{GMTW} and \cite[Proposition 4.2.4]{MW}, and rely on the surjectivety criteria. 

Let $X$ be a manifold without boundary and $A$ a closed subset in $X$. Let $s=(f,a)$ be a germ at $A$ of $\beta E^{\pitchfork}_{\J}(X)$; i.e., there is a locally finite cover $\{U_j\}$ of $A$ by open sets in $X$, a map $f$ from the union $U=\cup U_j$ to $\Omega^{\infty}_{\J}\MO$ and a broken function $a$ over $U$. 

Suppose that $f$ extends over $X$. In other words, the extended map $f$ represents an element in $E_{\J}[X, A; \alpha(s)]$. We need to show that it lifts to an element in $\beta E^{\pitchfork}_{\J}[X, A; s]$, i.e., the restriction $a|V$ to a sufficiently small open neighborhood $V$ of $A$ extends to a broken function over $X$. Let $O$ be an open set in $X$ such that the closures $\bar{O}$ and $\bar{V}$ are disjoint and $O\cup U=X$. By Lemma~\ref{ln:1}, we may perturb $f$ relative to $\bar{V}$ so that a representative of the adjoint of $f|O$ is a map 
\[  
     \tilde{f}: O\times \R\times S^{n+d-1}\longrightarrow {\mo}_n 
\]
smooth near $\tilde{f}^{-1}(\BO_n)$ and transversal to $\BO_n$. There is a continuous function $q$ over $\bar{O}$ that is greater than $a$ over $O\cap U$.  Now, for each $x$ in $\bar{O}$ we may choose a small neighborhood $U_x$ about $x$ and a function $a_x>q$ on $U_x$ such that the restriction of a representative of $f(y)$ to $\{a_x(y)\}\times S^{n+d-1}$ is transversal to $\BO_n$ for each $y\in U_x$. Thus we can choose a locally finite open covering of $\bar O$ together with a broken function $\{a_x\}$ over a neighborhood of $\bar O$. This defines a lift of $f$ to an element $(f, a)$ in $\beta E^{\pitchfork}_{\J}[X, A; s]$.

\end{proof}

\subsection{Weak homotopy equivalence 
$|\beta E^{\pitchfork}_{\J}| \simeq B|E^{\pitchfork}_{\J}|$}

Let $\mathcal{F}$ be a category valued sheaf on $\mathcal{X}$. It determines a set valued sheaf $N_0\mathcal{F}$ on $\mathcal{X}$ that takes a manifold $X$ onto the set of objects of the category $\mathcal{F}(X)$. Similarly $\mathcal{F}$ determines a set valued sheaf $N_1\mathcal{F}$ that assigns the set of morphisms of $\mathcal{F}(X)$ to a manifold $X\in \mathcal{X}$. In fact $\mathcal{F}$ determines a topological category denoted $|\mathcal{F}|$. Its space of objects is given by the classifying space $|N_0\mathcal{F}|$ of the set valued sheaf $N_0\mathcal{F}$, and the space of morphisms is given by $|N_1\mathcal{F}|$. By \cite[Theorem 4.1.2]{MW}, there is a homotopy equivalence 
\[
      |\beta \mathcal{F}| \simeq B|\mathcal{F}|. 
\]

We apply the Madsen-Weiss theorem in the case of the category valued sheaf $E^{\pitchfork}_{\J}$ to obtain a homotopy equivalence
\[
       |\beta E^{\pitchfork}_{\J}| \simeq B|E^{\pitchfork}_{\J}|.
\]

\subsection{Category valued sheaf $C_{\J}$}

Each topological category $\mathcal{C}$ gives rise to a category valued sheaf $C$ on $\mathcal{X}$ that takes a manifold $X$ onto the topological category $C(X)$ with space of objects given by 
the space
\[
    \Ob C(X)= C^{0}(X, \Ob \mathcal{C}),
\] 
of continuous maps into the topological space of objects of $\mathcal{C}$. The space of morphisms of $C(X)$ is the space  
\[
    \Mor C(X)= C^{0}(X, \Mor \mathcal{C}).
\] 
In particular the singular cobordism category $\mathcal{C}_{\J}$ gives rise to a category valued sheaf, which we denote by $C_{\J}$.

\subsection{Weak homotopy equivalence $B|E^{\pitchfork}_{\J}|\leftarrow B|C_{\J}|$}

Given a point $s$ in the space $\Ob C_{\J}(X)$, the composition
\[
   X\stackrel{s}\longrightarrow \Ob \mathcal{C}_{\J}\stackrel{(\ref{eq:10.1})}\longrightarrow \Omega^{\infty-1}_{\J}\MO\times \R,
\] 
determines an object $(f, a)$ in $E^{\pitchfork}_{\J}(X)$; the first component of the composition is the map $f$, while the second component is the function $a$. Similarly, given a morphism $s$, the composition 
\[
   X\stackrel{s}\longrightarrow \Mor \mathcal{C}_{\J}\stackrel{(\ref{eq:10.3})}\longrightarrow \Omega^{\infty-1}_{\J}\MO\times \R\times \R \times \R^2_{+}\longrightarrow \Omega^{\infty-1}_{\J}\MO\times \R^2_+,
\] 
where the last map is the obvious projection, determines a morphism in $E^{\pitchfork}_{\J}(X)$. In particular, we obtain a map of category valued sheaves
\[
  \gamma\colon   C_{\J}\longrightarrow E^{\pitchfork}_{\J}.
\]
As in \cite{GMTW}, one can show that the following lemma holds true. 

\begin{lemma}\label{l:10.2} There is a weak homotopy equivalence $B|C_{\J}|\to B|E^{\pitchfork}_{\J}|$. 
\end{lemma}
\begin{proof} For a category valued sheaf $\mathcal{C}$, the classifying space $B|\mathcal{C}|$ is the realization of the simplicial set
\[
       [k]\mapsto N_k|\mathcal{C}|,
\]
where $N_k$ is the $k$-th nerve of the category $|\mathcal{C}|$. In view of $N_k|\mathcal{C}|=|N_k\mathcal{C}|$, the above simplicial set is equivalent to 
\[
     [k]\mapsto |N_k\mathcal{C}|. 
\]
Thus, to prove Lemma~\ref{l:10.2} it suffices to show that $\gamma$ induces a weak equivalence of set valued sheaves $N_kC_{\J}$ and $N_kE^{\pitchfork}_{\J}$ for each $k\ge 0$. The proof of the latter is very similar to the proofs of \cite[Proposition 4.3]{GMTW} and \cite[Proposition 4.4]{GMTW}. Namely, we will make use of the surjectivety criteria. 

Let us show that for each manifold $X$, the map $\gamma$ induces a surjective map $N_kC_{\J}[X]\to N_kE^{\pitchfork}_{\J}[X]$ of concordance classes. An element $(f, \bar{a})$ in $N_kE^{\pitchfork}_{\J}[X]$ is represented by a map $f$ of $X$ to $\Omega^{\infty-1}_{\J}$ and a sequence of functions $a_0\le a_1 \le ... \le a_k$ on $X$. To simplify the argument, we may modify $(f,\bar{a})$ by concordance so that each function $a_i$ is constant. More importantly, we may further modify $(f, \bar{a})$ so that the map $f$ is represented by a map $\tilde f(t, s)=f_t(s)$,
\[
    f_t: S^{n+d-1}\longrightarrow {\mo}_n,   \qquad  t\in \R, 
\]
with $f_t(s)=f_{a_i}(s)$ for each $t$ in an $\varepsilon$-neighborhood of the value $a_i(X)$ for some $\varepsilon>0$. Finally we choose a diffeomorphism $\alpha\colon (a_0(X)-\varepsilon/2, a_k(X)+\varepsilon/2)\to \R$ which extends the inclusion of $(a_0(X)+\varepsilon/2, a_k(X)-\varepsilon/2)$. We can modify $(f, \bar{a})$ by concordance so that the representative $f_t(s)$ of $f$ is replaced by $f_{\alpha^{-1}(t)}(s)$. The obtained element represents an element in $N_kC_{\J}[X]$. 

The relative case is similar.   

\end{proof}

\subsection{Weak homotopy equivalence $B|C_{\J}|\to B\mathcal{C}_{\J}$} The proof is almost identical to the proof of \cite[Proposition 2.9]{GMTW} and uses the fact that the singular simplicial set of a topological space $X$ is weakly homotopy equivalent to $X$. 

\begin{lemma} There is a weak homotopy equivalence $B|C_{\J}|\to B\mathcal{C}_{\J}$. 
\end{lemma}
\begin{proof} The $k$-th nerve $N_k|C_{\J}|=|N_kC_{\J}|$ of the category $|C_{\J}|$ is the realization of the simplicial set
\[
       [l] \longrightarrow N_kC_{\J}(\Delta^l_e)=C^{0}(\Delta^l_e, N_k\mathcal{C}_{\J}),
\] 
and hence, for each $k$, the map $N_k|C_{\J}|\to N_k\mathcal{C}_{\J}$
is a weak homotopy equivalence, which implies the statement of the lemma. 
\end{proof}

This completes the proof of Theorem~\ref{th:1.4}. 

\begin{remark} An oriented version of Theorem~\ref{th:1.4} holds true as well. Namely, let $\J$ be an open stable set of singularities. Then there is a weak homotopy equivalence $B\mathcal{C}^+_{\J}\simeq \Omega^{\infty-1}\mathbf{B}^+_{\J}$. 
\end{remark}

\addcontentsline{toc}{part}{Decomposition $\Omega B\mathcal{C}_{\J}=\sqcup \BDiff M_{\alpha}$}

\section{Relation of singular cobordism categories to the b-principle}\label{s:last}
\subsection{Cobordism groups of singular maps}
We say that two proper $\J$-maps $f_i:M_i\to N$ of smooth manifolds, with
$i=0,1$,  are {\it $\J$-cobordant}, if there is a proper $\J$-map $F: W\to
N\times [0,1]$ of a manifold $W$ with boundary $M_0\sqcup M_1$
such that   
\[
     F(M_i)\subset N\times \{i\}   \qquad \textrm{for } i=0,1,
\] 
and the restrictions of $F$ to collar neighborhoods of $M_0$ and $M_1$ can be identified with the disjoint union of suspensions of $f_0$ and $f_1$. The set of $\J$-cobordant
classes of $\J$-maps of closed manifolds into a closed manifold $N$ gives rise to an abelian semigroup;  
the sum of elements $[f_0]$ and $[f_1]$ represented by maps $f_i: M_i\to N$, with $i=0,1$, is an element represented by the composition
\[
     M_0\sqcup M_1 \xrightarrow{f_0\sqcup f_1} N\sqcup N \longrightarrow N,
\]
where the last map is the identity map on each copy of $N$. In particular, the zero element is given by the map of an empty set. The semigroup of $\J$-cobordism classes of $\J$-maps into a manifold $N$ is turned into a group  $B^d(N; \J)$ by means of the Grothendieck construction; every element of $B^d(N; \J)$ is of the form $[f]-[g]$, where $[f]$ and $[g]$ are represented by $\J$-maps to $N$.

\emph{Oriented cobordism groups} of $\J$-maps can be defined similarly. 

A priori cobordism groups of $\J$-maps do not form a generalized cohomology theory since, for example, cobordism groups of $\J$-maps are not defined for topological spaces. 
There is, however, a 
counterpart of $B^d(N; \J)$ that can be used to compute $B^d(N; \J)$ in the same way as singular cohomology groups $H^n(N; \R)$, defined for every topological space, can be used to compute De Rham cohomology groups $H^n_{\mathrm{DR}}(N)$, defined only for smooth manifolds (see Theorem~\ref{th:0.2} below). 

Singular cobordism categories are related to cobordism groups of maps with prescribed singularities via a bordism version of the h-principle, which we state now.  We will omit here the precise definition of $\K$-invariance of sets of singularity types. Informally, a set $\J$ is {\it $\mathcal{K}$-invariant} if it has sufficiently many symmetries (e.g., see \cite{PW}). 

\begin{theorem}\label{th:0.2} {\rm{(}Sadykov, \cite{Sa}\rm{)}} 
Let $\J$ be an open, $\K$-invariant set of singularity types of map germs of dimension $d\ge 0$. Suppose that $\J$ contains all Morse singularity types. Then
\[  
      B^d(N; \J) \cong [N, \Omega^{\infty}\mathbf{B}_{\J}]
\]
for every closed manifold $N$ of dimension $\dim N>1$.  
\end{theorem}

\begin{remark} In the case $d<0$, which we do not consider in the current paper, the statement of Theorem~\ref{th:0.2} is true without the assumption that $\J$ contains all Morse singularity types \cite{An}, \cite{Sa},  \cite{Sz} (see also \cite{We} and \cite{El}). On the other hand, in the case $d\ge 0$ the statement of Theorem~\ref{th:0.2} is not true without the assumption that $\J$ contains all Morse singularity types for example for $\J=\emptyset$. 
\end{remark}

It is known~\cite{Mat} that a stable set $\J$ of singularity types is $\K$-invariant. 

Now we can identify the space $\Omega^{\infty}_{\J}\MO$ with the loop space $\Omega\Omega^{\infty-1}_{\J}\MO$. 

\begin{theorem}\label{th:10.2}
Let $\J$ be an open stable set of singularity types of map germs of dimension $d\ge 0$. Suppose that $\J$ contains Morse singularity types. Then 
\[
       \Omega^{\infty}_{\J}\MO\simeq \Omega B\mathcal{C}_{\J}.
\]
\end{theorem}
\begin{proof} The statement follows from the composition of weak homotopy equivalences
\[
     \Omega B\mathcal{C}_{\J}\longrightarrow \Omega^{\infty}\mathbf{B}_{\J} \longrightarrow \Omega^{\infty}_{\J}\MO
\]
where the first map is obtained from the map of Theorem~\ref{th:1.4} and the second map is the equivalence of two classifying spaces of cobordism groups of maps with prescribed singularities. More precisely, in order to construct the second weak homotopy equivalence, we again apply the argument in the proof of Theorem~\ref{ln:2}, but this time in the argument we refer not to Theorem~\ref{th:i1} but to Theorem~\ref{th:10.2}. We note that here we can not directly apply Theorem~\ref{ln:2} because a priori the space $\Omega^{\infty}_{\J}\MO$ and the loop space of $\Omega^{\infty-1}_{\J}\MO$ are not weakly homotopy equivalent. 
\end{proof}

\begin{remark} An argument similar to that in the proof of Theorem~\ref{th:10.2} shows that under the assumptions of Theorem~\ref{th:10.2}, there is a weak homotopy equivalence $\Omega^{\infty}_{\J}\MSO\simeq \Omega B\mathcal{C}^+_{\J}$.
\end{remark}

\section{Parametrized singular manifolds}\label{s:12}

Let $\J$ be a stable set of singularity types of map germs of dimension $d\ge 0$. Let $n$ be a fixed positive integer. An \emph{(embedded) parametrized singular manifold} with $\J$-singularities is represented by a smooth proper map 
\[
    \mathbf{f}\colon W\to D\times S^{n +d}
\]
where the first component $f$ is a proper $\J$-map into an open disc $D$ and the second component $j$ is an embedding whose image is disjoint from a compact neighborhood of the distinguished point in $S^{n+d}$. 
For $i=1,2$, two smooth maps $\mathbf{f}_i\colon W\to D\times S^{n+d}$ represent the same  parametrized singular manifold if the sets $f^{-1}_1(0)$ and $f_2^{-1}(0)$ are given by the same set $F\subset W$, and there is a neighborhood $U$ of $F$ in $W$ such that $\mathbf{f}_1|U=\mathbf{f}_2|U$. We say that two parametrized singular manifolds are of the \emph{same diffeomorphism type} if their representatives $\mathbf{f}_1$ and $\mathbf{f}_2$ can be chosen so that there are diffeomorphisms $\varphi$ and $\psi$ of pairs of spaces that fit a commutative diagram of maps
\[
\begin{CD}
(U_1, f_1^{-1}(0))@>\varphi >> (U_2, f_2^{-1}(0)) \\
@Vf_1 VV @Vf_2 VV \\
(V_1, 0)@>\psi>> (V_2, 0),
\end{CD}
\]
where $V_i$ is an open neighborhood of $0$ in $D$, and $U_i=f^{-1}_i(V_i)$ is an open subset in $W_i$. To simplify the notation we will use the same symbols for parametrized singular manifolds and their representatives. 
The pair $(\varphi, \psi)$ above represents an \emph{equivalence germ} between $\mathbf{f}_1$ and $\mathbf{f}_2$. Two pairs $(\varphi, \psi)$ and $(\varphi',\psi')$ represent the same germ if there are open neighborhoods $V_1$ and $U_1$ as above such that $\varphi=\varphi'$ over $U_1$ and $\psi=\psi'$ over $V_1$. If $\mathbf{f}_1=\mathbf{f}_2=\mathbf{f}$, then the set of equivalence germs is a group, denoted $\BDiff \mathbf{f}$, with operation given by taking the composition of equivalence germs.

\subsection{Stabilizations} In view of the inclusion of $D\times S^{n+d}$ into $D\times S^{n+d+1}$, every  parametrized singular manifold with $\J$-singularities  
\[
   \mathbf{f}\colon W\longrightarrow D\times S^{n+d}
\]
gives rise to a  parametrized singular manifold with $\J$-singularities
\[
    \mathbf{f}'\colon W\longrightarrow D\times S^{n+d+1}. 
\]
We take the colimit of the sets $\mathop{\mathrm{Set}}(\J, n)$ of parametrized singular manifolds  into $D\times S^{n+d}$ with respect to inclusions
\[
    \mathop{\mathrm{Set}}(\J, n)\longrightarrow \mathop{\mathrm{Set}}(\J, n+1),
\]
\[
    \mathbf{f}\mapsto \mathbf{f}'.
\]
Thus, by taking the colimit, we identify a parametrized singular manifold with a map 
\[  
   \mathbf{f}\colon W\longrightarrow D\times S^{\infty+d}
\]
whose image is in a manifold $D\times S^{n+d}$ for some finite $n$. 

We say that two parametrized singular manifolds $\mathbf{f_i}\colon W_i\to D_i\times S^{\infty+d}$ with $i=1,2$ are of the \emph{same stable diffeomorphism type} if there are non-negative integers $n_1$ and $n_2$ and diffeomorphisms $\varphi$ and $\psi$ that fit a commutative diagram of maps
\[
\begin{CD}
(U_1\times \R^{n_2}, f_1^{-1}(0)\times \{0\})@>\varphi >> (U_2\times \R^{n_1}, f_2^{-1}(0)\times \{0\}) \\
@Vf_1\times \id_{\R^{n_2}} VV @Vf_2\times \id_{\R^{n_1}} VV \\
(V_1\times \R^{n_2}, \{0\}\times \{0\})@>\psi>> (V_2\times \R^{n_1}, \{0\}\times\{0\}),
\end{CD}
\]
where $V_i$ is an open neighborhood of $0$ in $D_i$ and $U_i=f^{-1}_i(V_i)$. 

A priori it is possible that two parametrized singular manifolds $\mathbf{f_1}$ and $\mathbf{f_2}$ are of different diffeomorphism types, while the parametrized singular manifolds $\mathbf{f_1}\times \id_{\R^n}$ and $\mathbf{f_2}\times \id_{\R^n}$ are of the same diffeomorphism type for some $n$. To avoid this, we will consider only \emph{stable singular manifolds}, which we define next. In terms of representatives an \emph{unfolding} of a parametrized singular manifold 
\[
   \mathbf{f}\colon W\longrightarrow D\times S^{n+d}
\]
is a smooth family $\mathbf{F}$ of parametrized singular manifolds  
\[
      \mathbf{F}_t\colon W\longrightarrow D\times S^{n+d}
\]
where $t\in D^k$ and $\mathbf{F}_0=\mathbf{f}_0$. We note that $\mathbf{F}$ is itself a parametrized singular manifold,  
\[
      \mathbf{F}\colon D^k \times W\longrightarrow (D^k\times D)\times S^{n+d}
\]
We say that a parametrized singular manifold $\mathbf{f}$ is \emph{stable} if for each unfolding $\mathbf{F}$ of $\mathbf{f}$ there is a smooth family of equivalence maps $(\varphi_t, \psi_t)$ between the unfolding $\mathbf{F}$ and the trivial unfolding $\mathbf{f}\times \id_{D^k}$. 

For each stable diffeomorphism type $M_{\alpha}$ of  parametrized singular manifolds there is a parametrized singular manifold given by a map $\mathbf{f}$ of a manifold of minimal dimension. We define the group $\Diff M_{\alpha}$ to be the group $\Diff\mathbf{f}$. Let $\mathop{\mathrm{EDiff}}M_{\alpha}$ denote the topological space of  parametrized singular manifolds of the same diffeomorphism type $\mathbf{f}$, where $\mathbf{f}$ is the minimal  parametrized singular manifold of the type $M_{\alpha}$; the space $\mathop{\mathrm{EDiff}}M_{\alpha}$ is endowed with the Whitney $C^{\infty}$-topology. It is weakly homotopy equivalent to a point. 

The group $\Diff M_{\alpha}$ acts on the space $\mathop{\mathrm{EDiff}}M_{\alpha}$. Namely, an equivalence germ $(\varphi, \psi)$ takes a parametrized singular manifold $\mathbf{f}$ onto $\psi\circ\mathbf{f}\circ\varphi^{-1}$. The quotient with respect to this action is the space denoted $\BDiff M_{\alpha}$. The projection
\[
   \mathop{\mathrm{EDiff}}M_{\alpha} \longrightarrow \BDiff M_{\alpha}
\]
is a fiber bundle with fiber $\Diff M_{\alpha}$. Consequently the space $\BDiff M_{\alpha}$ is the classifying space for families of  singular manifolds of type $M_{\alpha}$. 

\begin{example}
In the case where the fiber $M_{\alpha}$ is non-singular, the group $\Diff M_{\alpha}$ is the diffeomorphism group of the manifold $M_{\alpha}$, while the space $\BDiff M_{\alpha}$ is the classifying space of $\Diff M_{\alpha}$. 
\end{example}

For a stable diffeomorphism type $M_{\alpha}$, there are homomorphisms of groups 
\[
    \Diff M_{\alpha}\longrightarrow \Diff D,
\]
\[
   (\varphi, \psi)\longrightarrow \psi,
\]
where $\Diff D$ is the group of diffeomorphism germs of $\Diff D$ at $0$, and
\[
    \Diff M_{\alpha}\longrightarrow \Diff W,
\] 
\[
   (\varphi, \psi)\longrightarrow \varphi,
\]
where $W$ is the source manifold of a minimal representative $\mathbf{f}$ of $M_{\alpha}$. In particular, the group $\Diff M_{\alpha}$ acts on $W$ and $D$. We put
\[
   \mathcal{E}M_{\alpha}\colon=\mathop{\mathrm{EDiff}}M_{\alpha}\times_{\Diff M_{\alpha}} W, 
\]
\[
   \mathcal{B}M_{\alpha}\colon=\mathop{\mathrm{EDiff}}M_{\alpha}\times_{\Diff M_{\alpha}} D.
\]
Then the \emph{universal family} of stable singular manifolds of type $M_{\alpha}$ is defined to be the family

\[
\xymatrix{
  &  \mathcal{E}M_{\alpha} \ar[rr] \ar[dr] &     &  \mathcal{B}M_{\alpha} \ar[dl] & \\
  &                                 & \BDiff M_{\alpha}   &                     &
}
\]


\subsection{From singular manifolds to free singular manifolds}

Every  paramet\-rized singular manifold $\mathbf{f}$ determines a compact subset $F=f^{-1}(0)$ embedded into $S^{n+d}$ by means of $j$. The compact subset $F\subset S^{n+d}$ alone does not determine the diffeomorphism type of the singular manifold $\mathbf{f}$. However, in this subsection we will show that the diffeomorphism type is determined if $F$ is equipped with an additional structure (see Theorem~\ref{th:6.2}).  

Let $Y_i$ be a compact subset of a smooth manifold $M_i$ for $i=1,2$. We recall that a continuous map $Y_1\to Y_2$ is said to be a \emph{diffeomorphism} if it extends to a diffeomorphism of open neighborhoods of $Y_1$ in $M_1$ and $Y_2$ in $M_2$. Similarly a vector bundle over a compact subset $Y$ of a manifold $M$ is called \emph{smooth} if it can be extended to a smooth vector bundle over an open neighborhood of $Y$ in $M$.  

\begin{remark}[Motivation for Definition~\ref{d:3}]
Let $\mathbf{f}$ be a parametrized singular manifold. Then, for each point $x\in F$, it determines a map germ from a neighborhood of $x$ in $W$ to $D$. Let $\mathcal{O}(F)$ denote a neighborhood of the zero section in the tangent bundle of $W$ restricted to $F$. For each point $x\in F$ we may use a Riemannian metric on $W$ to identify the fiber $\mathcal{O}_xF$ of $\mathcal{O}(F)$ over $x$ with a neighborhood of $x$ in $W$. Thus we deduce that every  parametrized singular manifold $\mathbf{f}$ determines a smooth family $\mathfrak{o}: \mathcal{O}(F)\to D$ of map germs  $\mathcal{O}_xF\to D$ parametrized by points $x$ in $F$. Similarly we observe that $\mathbf{f}$ determines a smooth family $\mathfrak{u}: \mathcal{U}\to \mathcal{O}(F)$ of proper embedding germs $\mathfrak{u}_x\colon \mathcal{U}_x\to\mathcal{O}_xF$, where $\mathcal{U}_x$ is an open neighborhood of $x$ in $F$. 
\end{remark}

\begin{definition}\label{d:3} Let $F$ be a closed subset in $S^{n+d}$. A \emph{free singular manifold} structure on $F$ of dimension $d$ is a triple $(\xi, \mathfrak{o}, \mathfrak{u})$ of 
\begin{itemize}
\item a smooth vector bundle $\xi$ over $F$ of dimension $m+d$, 
\item a smooth family $\mathfrak{o}$ of map germs 
\[
     \mathfrak{o}_x\colon (\mathcal{O}_xF, x) \longrightarrow (D, 0)
\]
parametrized by $x\in F$, where $\mathcal{O}_xF$ is a neighborhood of $x$ in the fiber $\xi_x$ and $D$  is a disc of dimension $m$; the family $\mathfrak{o}$ is represented by a smooth map from an open neighborhood $\mathcal{O}(F)$ of the zero section of $\xi$ to $D$, and
\item a smooth family $\mathfrak{u}$ of proper embedding germs 
\[
    \mathfrak{u}_x\colon (\mathcal{U}_x, x)\longrightarrow (\mathcal{O}_xF, x)
\]
parametrized by $x\in F$, where $\mathcal{U}_x$ is an open neighborhood of $x$ in $F$ such that $\mathfrak{o}^{-1}_x(0)=\mathfrak{u}_x(\mathcal{U}_x)$. 
\end{itemize} 
\end{definition}

\begin{theorem}\label{th:6.2} Every subset $F\subset S^{n+d}$ with a free singular manifold structure $(\xi, \mathfrak{o}, \mathfrak{u})$ determines a parametrized singular manifold $\mathbf{f}=(f, j)$, unique up to a choice of the $j$-component.  
\end{theorem}
\begin{proof} To simplify the argument we assume that $F$ has only one singular point $s$. In general the argument is similar. For each point $x\ne s$ in $F$ we fix a small neighborhood $U_x$ of $x$ in $F$ so that $s\notin U_x$. Next for each point $x$ in $F$ we fix an open neighborhood $O_x(F)$ in $\xi_x$ so that for $x\ne s$ the set $O_xF$ consists of geodesic discs $\mathcal{D}_t$, where $t\in U_x$ and $\mathcal{D}_t$ is given by the union of geodesic rays of short length emitted from $\mathfrak{u}_x(t)$ in the direction normal to $d\mathfrak{u}_x(T_tF)$. There is a smooth map 
\[
     (y_x, z_x)\colon (O_xF, x)\longrightarrow (F\times D, \{x\}\times \{0\})
\] 
where $y_x$ maps each geodesic disc $\mathcal{D}_t$ onto $t\in F$ and $z_x$ coincides with $\mathfrak{o}_x|O_xF$. We may assume that for $x$ close to $s$ the geodesic discs $\mathcal{D}_t$ in $O_xF$ map to the corresponding geodesic discs $\mathcal{D}_t$ in $O_sF$. Finally, we define the manifold $W$ to be the quotient
\[
      (\bigsqcup_{x\ne s} O_xF)\sqcup O_sF/\sim
\]  
where a point $a$ in the fiber $O_xF$, with $x\ne s$ is equivalent to a point $b$ in the fiber $O_{x'}F$, with $x'\ne s$, if and only if 
\[
 y_x(a) =  y_{x'}(b)
\]
and
\[
    z_x(a)=z_{x'}(b). 
\]
Points in $O_sF$ are identified with points in $O_xF$ if there is a smooth map $i\colon O_xF\to O_sF$ that restricts to an inclusion 
\[
 \mathfrak{u}_x(U_x)\stackrel{(\mathfrak{u}_x|U_x)^{-1}}\longrightarrow U_x\stackrel{\subset}\longrightarrow U_s \stackrel{\mathfrak{u}_x}\longrightarrow \mathfrak{u}_x(U_s)
\]
and maps every geodesic disc $\mathcal{D}_t$ in $\mathcal{O}'_xF$ normal to $\mathfrak{u}_x(\mathcal{U}_x)$ to a geodesic disc in $\mathcal{O}_sF$ normal to $\mathfrak{u}_x(\mathcal{U}_s)$ so that $\mathfrak{o}_x|\mathcal{D}_t=\mathfrak{o}_s\circ i|\mathcal{D}_t$.

The smooth structure on $W$ is determined by the cover of $W$ by open sets $\mathcal{O}'_xF$ and $\mathcal{O}_sF$. The map $f: W\to D$ is defined by the family $\mathfrak{o}$. On the other hand the embedding $j: W\to S^{n'+d}$ can be chosen to be an arbitrary extension of the smooth embedding $F\to S^{n+d}\subset S^{n'+d}$ provided that $n'$ is sufficiently big.    
\end{proof}

We say that two free singular manifolds are of the same diffeomorphism type if they determine parametrized fiber germs of the same diffeomorphism type. 

A \emph{family of free singular manifolds} parametrized by a manifold $N$ is a subset $N\times F\subset N\times S^{n+d}$ together with a $4$-tuple $(\xi, \gamma, \mathfrak{o}, \mathfrak{u})$ of
\begin{itemize}
\item a continuous vector bundle $\xi$ over $N\times F$ of dimension $m+\dim N$,
\item a continuous vector bundle $\gamma$ over $N$, 
\item a continuous family of map germs
\[
   \mathfrak{o}_{y\times x}\colon (\mathcal{O}_{y\times x}N\times F, \{y\}\times \{x\})\longrightarrow (\gamma_y, 0)
\]
where $\mathcal{O}_{y\times x}N\times F$ is an open neighborhood of the origin in the vector space $\xi_{y\times x}$, 
\item a continuous family of proper embedding germs 
\[
   \mathfrak{u}_x\colon (\mathcal{U}_{x,y}, \{y\}\times\{x\})\longrightarrow (\mathcal{O}_{y\times x}N\times F, \{y\}\times \{x\})
\]
where $\mathcal{U}_{y\times x}$ is a neighborhood of the point $y\times x$ in $\{y\}\times F$
\end{itemize}
such that 
\begin{itemize}
\item for each $x\in N$ the set $\{x\}\times F$ is embedded into $\{x\}\times S^{n+d}$, and
\item for each $y\in N$ the $4$-tuple restricts to a free singular manifold structure on $\{y\}\times F$, and
\item free singular manifold structures on $\{y\}\times F$ for all $y\in N$ are of the same diffeomorphism type.	 
\end{itemize}

\begin{theorem}\label{th:6.3} Every family of free singular manifolds determines a family of parametrized singular manifold, unique up to a continuous deformation of families of parametrized singular manifolds.  
\end{theorem}
\begin{proof} The proof follows from the construction in the proof of Theorem~\ref{th:6.2}. 
\end{proof}

\section{Decomposition of $\Omega B\mathcal{C}_{\J}$}\label{s:13}

Let us note that not each point in the space $\Omega^{\infty}_{\J}\MO$ is in the image of Pontrjagin-Thom maps. Nevertheless, each loop in this space determines not only a subset $F\subset S^{\infty}$ but also a structure of a free singular manifold on $F$.  

We recall that there are two vector bundles over the Grassmannian manifold $G_n(\R^{n+d})=G_d(\R^{n+d})$, one is of dimension $n$ with total space $E_{n}(\R^{n+d})$ and the other is of dimension $d$ with total space $E_{d}(\R^{n+d})$. Each point in $\Omega^{\infty}_{\J}\MO$ corresponds to the loop
\[
    f\colon     S^{n+d}\longrightarrow TE_n(\R^{n+d})
\]
smooth near $F=f^{-1}(G_n(\R^{n+d}))$. It determines the free singular manifold structure $(\xi, \mathfrak{o}, \mathfrak{u})$ on the set with $F$. Namely, here 
\begin{itemize}  
\item $\xi$ is a vector bundle over $F$ of dimension $n+2d$; it is given by the Whitney sum $\gamma^d\oplus \varepsilon^{n+d}$, where $\gamma^d$ stands for the pullback bundle $(f|F)^*E_d(\R^n)$ and $\varepsilon$ is the trivial bundle over $F$, 
\item $\mathcal{O}(F)$ is a small open neighborhood of the zero section in $\xi$; we use the Riemannian metrics to identify the fiber $\mathcal{O}_xF$ of $\mathcal{O}(F)$ over $x$ with the product of a neighborhood of the origin in the vector space $\gamma^d_x$ and a neighborhood of $x$ in $S^{n+d}$,
\item $D$ is a disc of dimension $n+d$,  
\item $\mathfrak{o}_x\colon \mathcal{O}_xF\to D$ is a map germ given by 
 \begin{eqnarray}
      (\gamma^d_x\times S^{n+d}, \{0\}\times \{x\})&\stackrel{\id \times f}\longrightarrow &(\gamma^d_x\times {\mo}_n, \{0\}\times \{f(x)\}) \nonumber \\
       &\longrightarrow & (\gamma^d_x\oplus \gamma^n_x, \{0\}\times \{0\}) \nonumber \\
       & \longrightarrow & (D, 0) \nonumber,
 \end{eqnarray}
where the second map is given by the product of the identity map on $\gamma^d_x$ and a projection onto the fiber $E_n(\R^{n+d})|f(x)$ which we identify with its pullback $\gamma^n_x$ with respect to $f$; and 
 \item $\mathfrak{u}_x$ is the inclusion of a neighborhood of $x$ in $F$ to 
 \[
 T_xS^{n+d}= \{0\}\times T_xS^{n+d}\subset \gamma^d_x \oplus \varepsilon^{n+d}_x.
 \]  
\end{itemize}

Let $M_{\alpha}$ be a stable singular manifold of maps of dimension $d$. Let $bM_{\alpha}$ denote the subset of $\Omega B\mathcal{C}_{\J}\simeq \Omega^{\infty}_{\J}\MO$ that consists of loops corresponding to free singular manifolds of stable singular manifold type $M_{\alpha}$. Since every point in $\Omega^{\infty}_{\J}\MO$ corresponds to a free singular manifold,  
\[
     \Omega^{\infty}_{\J}\MO=\bigsqcup_{M_{\alpha}} bM_{\alpha}.
\]
We observe that the above construction has a parametric generalization. Namely, every map $S\to bM_{\alpha}$ determines a family of singular manifolds of stable diffeomorphism type $M_{\alpha}$.

\begin{theorem}\label{th:7.1} The subspace $bM_{\alpha}$ is weakly homotopy equivalent to $\BDiff M_{\alpha}$. 
\end{theorem}
\begin{remark} In the case of maps of negative dimension $d$, a fiber consists at most of finitely many points. In this case (which we do not consider in the current paper) a singular manifold is called a \emph{multigerm} and Theorem~\ref{th:7.1} follows from the paper \cite{RS} by Rim\'anyi and Sz\H{u}cs where a completely different method is applied. 
\end{remark}

\begin{proof} To begin with we use a parametric version of the Pontrjagin-Thom construction to the universal family of maps $\mathcal{E}M_{\alpha}\to \mathcal{B}M_{\alpha}$ over $\BDiff M_{\alpha}$ to obtain a map 
\begin{equation}\label{eq:7.2}
     \mathcal{B}M_{\alpha}\longrightarrow \Omega^{\infty}_{\J}\MO.
\end{equation}
We note that the zero section of the disc bundle $\mathcal{B}M_{\alpha}$ is isomorphic to $\BDiff M_{\alpha}$ and its image under the map (\ref{eq:7.2}) is in $bM_{\alpha}$. Thus we obtain a map
\begin{equation}\label{eq:7.1}
    \BDiff M_{\alpha}\longrightarrow bM_{\alpha}. 
\end{equation}
Let us show that this map induces a surjective homomorphism in homotopy groups. Let $S\to bM_{\alpha}$ be a continuous map. By the construction before Theorem~\ref{th:7.1}, we obtain a family of free singular manifolds parametrized by $S$. Such a family leads to a continuous map $S\to \BDiff M_{\alpha}$ which is a lift of the map (\ref{eq:7.1}). Thus the induced map in homotopy groups is surjective. A similar argument shows that the induced map is injective. 
\end{proof}

Theorem~\ref{th:7.1} immediately implies Theorem~\ref{th:main}.

\begin{remark} An oriented version of the decomposition of $\Omega B\mathcal{C}^+_{\J}$ holds true, namely as a set the space $\Omega B\mathcal{C}^+_{\J}$ is given by $\cup \BDiff^+ M_{\alpha}$ where $M_{\alpha}$ is an oriented singular manifold of dimension $d$.  
\end{remark}

\end{document}